\documentclass[12pt,leqno]{article}
\usepackage{amsmath,amsthm,amscd,amssymb}
\usepackage{latexsym}
\numberwithin{equation}{section}

\newtheorem{thm}{Theorem}[section]
\newtheorem{lemma}{Lemma}[section]  \theoremstyle{remark}
\newtheorem*{rem}{Remark}
\newtheorem*{rems}{Remarks}

\def\beq{\begin{equation}}
\def\endq{\end{equation}}
\def\RR{\mathbb{R}}
\def\SS{\mathbb{S}}
\def\HH{\mathbb{H}}

\begin{document} 

\title
{A Sharp Bound for the Ratio of the First Two Dirichlet
Eigenvalues of a Domain in a Hemisphere of $\SS^n$}   

\author{Mark S. Ashbaugh\thanks{Partially
supported by National Science Foundation
(USA) grants DMS--9114162, INT--9123481, DMS--9500968, and 
DMS--9870156.}\\
Department of Mathematics\\ 
University of Missouri\\ 
Columbia, Missouri 65211--0001, USA\\
mark@math.missouri.edu 
\and  
Rafael D. Benguria\thanks{Partially
supported by FONDECYT (Chile) project number
196--0462 and a C\'atedra Presidencial en Ciencias (Chile).}\\
Departamento de F\'\i sica\\ 
P. Universidad Cat\'olica de Chile\\
Casilla 306, Santiago 22, CHILE\\ 
rbenguri@fis.puc.cl} 

\maketitle 

\begin{description}
\item[1991 Mathematics Subject Classification:] Primary 58G25,
Secondary 35P15, 49Rxx, 33C55.
\item[Keywords and phrases:] eigenvalues of the Laplacian, Dirichlet
problem for domains on spheres, Payne--P\'{o}lya--Weinberger
conjecture, Sperner's inequality, ratios of eigenvalues, isoperimetric
inequalities for eigenvalues. 
\item[Short title:  The PPW Conjecture in $\SS^n$] 
\end{description}   

\begin{abstract} 
For a domain $\Omega$ contained in a hemisphere of the 
$n$--dimensional sphere $\SS^n$ we prove the optimal result 
$\lambda_2/\lambda_1(\Omega) \le \lambda_2/\lambda_1(\Omega^{\star})$ 
for the ratio of its first two Dirichlet
eigenvalues where $\Omega^{\star}$, the symmetric rearrangement of 
$\Omega$ in $\SS^n$, is a geodesic ball in $\SS^n$ having the same 
$n$--volume as $\Omega$.  We also show that $\lambda_2/\lambda_1$ for 
geodesic balls of geodesic radius $\theta_1$ less than or equal to 
$\pi/2$ is an increasing function of $\theta_1$ which runs between the 
value $(j_{n/2,1}/j_{n/2-1,1})^2$ for $\theta_1=0$ (this is the Euclidean
value) and $2(n+1)/n$ for $\theta_1=\pi/2$.  Here $j_{\nu,k}$ denotes
the $k^{th}$ positive zero of the Bessel function $J_{\nu}(t)$.  This
result generalizes the Payne--P\'{o}lya--Weinberger conjecture, which
applies to bounded domains in Euclidean space and which we had proved
earlier.  Our method makes use of symmetric rearrangement of functions
and various technical properties of special functions. We also prove
that among all domains contained in a hemisphere of $\SS^n$ and
having a fixed value of $\lambda_1$ the one with the maximal value of
$\lambda_2$ is the geodesic ball of the appropriate radius.  This is a 
stronger, but slightly less accessible, isoperimetric result than
that for $\lambda_2/\lambda_1$.  Various other results for $\lambda_1$
and $\lambda_2$ of geodesic balls in $\SS^n$ are proved in the course
of our work.
\end{abstract}

\bigskip 

\section{Introduction}\label{s1}  

\bigskip

    With our earlier proof \cite{AB91}, \cite{AB92a}, \cite{AB92b} of the
Payne--P\'{o}lya--Weinberger conjecture \cite{PPW55}, \cite{PPW56} the bound
\begin{equation}
\lambda_2/\lambda_1(\Omega) \le \lambda_2/\lambda_1(\Omega^{\star})=
j_{n/2,1}^2/j_{n/2-1,1}^2   
\label{eq:I1}
\end{equation}
was established for the ratio of the first two eigenvalues of the Laplacian
$-\Delta$ on a bounded domain $\Omega \subset \RR^n$ with Dirichlet boundary
conditions imposed on $\partial \Omega$. Here $\Omega^{\star}$ represents 
the $n$--dimensional ball having the same volume as $\Omega$ (but, in fact, 
by scaling any ball will do) and $j_{\nu,k}$ represents the $k^{th}$ 
positive zero of the Bessel function $J_{\nu}(t)$ \cite{AS64}.  Equality 
obtains in (\ref{eq:I1}) if and only if $\Omega$ is a ball to begin with.
In this paper we prove the analog of this result for domains in a
hemisphere of $\SS^n$ ($-\Delta$ is now, of course, the Laplacian on $\SS^n$). 

It turns out to be better to view (\ref{eq:I1}) as
\begin{equation}
\lambda_2(\Omega) \le \lambda_2 (B_{\lambda_1})  
\label{eq:I2}
\end{equation}
where $B_{\lambda_1}$ represents the $n$--dimensional ball that has 
the value $\lambda_1(\Omega)$ as its first Dirichlet eigenvalue.  This 
is, in fact, the way our proof proceeded (\cite{AB91}, \cite{AB92a}, 
\cite{AB94b}, see also \cite{Chi3}).  Of course, the choice of ball 
here just involves choosing an appropriate radius and this choice is 
always unique since the first eigenvalue of a ball is a strictly 
decreasing function of its radius and goes from infinity to zero as 
the radius goes from zero to infinity.  In words, (\ref{eq:I2}) says 
that among all $n$--dimensional domains having the same first eigenvalue 
the $n$--dimensional ball has maximal second eigenvalue.  One might
compare this statement with the statement of the Faber--Krahn
inequality \cite{Fa23}, \cite{Kr25}, \cite{Kr26}:  among
all $n$--dimensional domains having the same volume the $n$--dimensional 
ball has minimal first (Dirichlet) eigenvalue.  Also of interest is the 
Szeg\H{o}--Weinberger inequality \cite{Sz54}, \cite{We56}:
among all $n$--dimensional domains having the same volume the 
$n$--dimensional ball has maximal first nonzero Neumann eigenvalue. 
Both of these other inequalities are relevant here; the first because 
it figures in our proof and the second because in many ways the proof 
of this result is analogous to (though considerably simpler than) our 
proof of the Payne--P\'{o}lya--Weinberger conjecture. 

What we do in this paper is transfer the strategy outlined in the last
paragraph over to bounded domains in $\SS^n$ contained in hemispheres. 
It turns out that, properly interpreted, everything that we have said 
above concerning domains in Euclidean space also holds in $\SS^n$.  Thus 
we prove that (\ref{eq:I2}) holds if $\Omega$ is a domain in a hemisphere 
in $\SS^n$ and $B_{\lambda_1}$ is the geodesic ball in $\SS^n$ having
$\lambda_1(\Omega)$ as its first eigenvalue.  There is a Faber--Krahn
result in $\SS^n$ \cite{Sp73} (see also \cite{FrHa76}), too, and
$\lambda_1$ of a geodesic ball is still strictly monotone decreasing
and goes to infinity at zero so a unique $B_{\lambda_1}$ which is a
hemisphere or less will exist.  For the Faber--Krahn and Szeg\H{o}--Weinberger
inequalities in $\SS^n$ one need only read ``volume'' as ``canonical
$n$--dimensional volume in $\SS^n$'' and ``ball'' as ``geodesic
ball''.  The generalization of the Szeg\H{o}--Weinberger result to
domains in hemispheres of $\SS^n$ is a recent result of ours \cite{AB95} 
(see also prior work of Chavel \cite{Cha1}).  Our proof here parallels the
one in \cite{AB95} as well as our proof of the
Payne--P\'{o}lya--Weinberger conjecture \cite{AB91}, \cite{AB92a},
\cite{AB92b}  (see especially our proof in \cite{AB92b}) for the
Euclidean case. 
   
The biggest difference between our work on $\lambda_2/\lambda_1$ for 
domains in $\SS^n$ versus those in $\RR^n$ revolves around the difference
between (\ref{eq:I1}) and (\ref{eq:I2}). 
In $\RR^n$ (\ref{eq:I1}) and (\ref{eq:I2}) are equivalent since
$\lambda_2/\lambda_1$ is the same for any ball, whatever its size.  This
follows from the fact that in that case the eigenvalues scale with the 
radius, in particular, for a ball of radius $R$ in $\RR^n$ {} $\lambda_1
=j^2_{n/2-1,1}/R^2$, $\lambda_2 = j^2_{n/2,1}/R^2$, and hence
$\lambda_2/\lambda_1 =j_{n/2,1}^2/j_{n/2-1,1}^2$ and is independent of $R$.
When one passes to $\SS^n$ this is no longer the case.  If we let $\theta_1$
denote the radius of our geodesic ball then the first and second eigenvalues
of that ball, which we denote by $\lambda_1(\theta_1)$ and
$\lambda_2(\theta_1)$, are in general more complicated functions of
$\theta_1$.  By domain monotonicity (see, for example,
\cite{Ba80}, \cite{Cha2}, \cite{CoHi}, \cite{Rayl}) these are, of
course, strictly decreasing functions but more precise knowledge of them
(or combinations thereof) requires considerable effort.  For
example, for the case of $\SS^n$, to pass from (\ref{eq:I2}) back to the 
first part of (\ref{eq:I1}) (i.e., $\lambda_2/\lambda_1(\Omega) \le 
\lambda_2/\lambda_1(\Omega^{\star})$, where $\Omega^{\star}$ is the
symmetric rearrangement of $\Omega$ on $\SS^n$) one needs to know that
$\lambda_2/\lambda_1$ for geodesic balls is an increasing function of
$\theta_1$.  We do this below in Section 3.
In general, most of what could almost be taken for granted for the
case of balls in $\RR^n$ expands to some problem about how
$\lambda_1$, $\lambda_2$, or a combination of the two behaves as a
function of $\theta_1$.  Another instance of this is that while for a
ball in $\RR^n$ it is easy to see that $\lambda_2$ corresponds to an
$\ell =1$ eigenfunction (i.e., an eigenfunction associated with an
$\ell =1$ spherical harmonic) in $\SS^n$ the analog of this must be
proved for all $\theta_1 \in (0,\pi/2]$.  The proof of this fact, while 
not difficult, is given in Section 3.  In fact, in Section 3 (see Lemma 
\ref{lem:3.1}) we prove the result for all $\theta_1 \in (0,\pi)$ 
(note that $m$ replaces $\ell$ there).  

For orientation we outline the elements of our proof here.  Aside from 
item 4 (which we have just discussed) these elements were also present 
in our proof in the Euclidean case. 

\begin{enumerate} 

\item 
Rayleigh--Ritz inequality for estimating $\lambda_2$:
$$
\lambda_2-\lambda_1 \le \frac{\int_{\Omega}{\vert \nabla P \vert}^2 u_1^2 
\, d\mu}{\int_{\Omega} P^2 u_1^2 \, d\mu}
$$
if $\int_{\Omega} P u_1^2 d\mu=0$ and $P \not\equiv 0$.  Here $u_1$ denotes 
the normalized eigenfunction for $\lambda_1$.  This inequality applies on 
manifolds as in Euclidean space; one only has to view $\vert \nabla P \vert$ 
as a norm with respect to the metric on the manifold (see, e.g., 
\cite{Cha2}, pp.~50--51).  Also $d\mu$ represents the intrinsic volume 
element for the manifold.  The inequality given here follows from the 
usual Rayleigh-Ritz inequality for $\lambda_2$ by taking $P u_1$ as trial 
function and integrating by parts appropriately.  

\item
 A Brouwer fixed point theorem argument that allows us to insure that
the condition $\int_{\Omega} Pu_1^2 d \mu=0$ is satisfied for $n$
specific choices of the function $P$.  In our recent paper \cite{AB95} we 
gave a version of this argument for $\SS^n$ which also applies here. 
In Section 2 of this paper we give an improved version of this
argument (using degree theory rather than the Brouwer fixed point
theorem) which yields a slightly stronger result.  The 
original argument of this type (to our knowledge) was given by
Weinberger in \cite{We56}.  For future reference, we note that all 
such results will be referred to as ``center of mass results''.  

\item
Rearrangement results for functions and domains.  These results are
measure theoretic in nature and easily extend to problems on manifolds. 
There is also an easily effected preliminary rearrangement in $\SS^n$ (see
also \cite{AB94b}, \cite{AB95}).  

\item
Properties of the eigenvalues and eigenfunctions of geodesic balls
in $\SS^n$.  These we establish in Section 3 below.  The details depend on
Legendre and associated Legendre functions though we manage to keep these
functions (also expressible in terms of hypergeometric functions) in the
background.  In $\RR^n$, of course, the special functions that occur are 
Bessel functions.  

\item
Monotonicity properties of certain special combinations of the
eigenfunctions for $\lambda_1$ and $\lambda_2$ for geodesic balls
contained in hemispheres.  These functions are special to the problem
of maximizing $\lambda_2/\lambda_1$ (or $\lambda_2$ subject to
$\lambda_1=const.$) and in a technical sense their properties proved
here are the most difficult part of our overall proof.  These
properties are proved in Section 4 below. 

\item
Chiti's comparison argument.  This is a specialized comparison result
which establishes a crossing property of the symmetric--decreasing
rearrangement of the eigenfunction $u_1$ vis--\`a--vis the first 
eigenfunction of the geodesic ball $B_{\lambda_1}$.  This result of 
Chiti \cite{Chi0}, \cite{Chi1}, \cite{Chi2}, \cite{Chi3} for the 
Euclidean case is based upon a rearrangement technique for partial
differential equations due to Talenti \cite{Ta76b} which in turn is 
based on the classical isoperimetric inequality in $\RR^n$ \cite{BuZa88}, 
\cite{Cha2}, \cite{Os78}.  All of these results generalize to $\SS^n$.  
We give these arguments in detail in Section 5. 
\end{enumerate}  

\bigskip

Our main results are summarized in the following theorems: 

\bigskip  
\begin{thm} \label{thm:I1}
Let $\Omega$ be contained in a hemisphere of $\SS^n$ and let
$B_{\lambda_1}$ denote the geodesic ball in $\SS^n$ having the same
value of $\lambda_1$ as $\Omega$ (i.e., $\lambda_1 (\Omega)
=\lambda_1(B_{\lambda_1})$).  Then 
$$
\lambda_2(\Omega) \leq \lambda_2(B_{\lambda_1})
$$
with equality if and only if $\Omega$ is itself a geodesic ball in $\SS^n$.
\end{thm}  

\bigskip  

\begin{thm} \label{thm:I2}
The first Dirichlet eigenvalue for a geodesic ball in
$\SS^n$ of geodesic radius $\theta_1$, $\lambda_1(\theta_1)$, is such that
$\theta_1^2\lambda_1(\theta_1)$ is a decreasing function of
$\theta_1$ for $0<\theta_1 \le \pi$.
\end{thm}  

\bigskip  

\begin{thm}  \label{thm:I3}
The quotient between the first two Dirichlet eigenvalues for
a geodesic ball in $\SS^n$, of geodesic radius $\theta_1$,
$\lambda_2/\lambda_1$, is an increasing function of $\theta_1$ for
$0<\theta_1 \le \pi/2$.
\end{thm}

Finally, we come to the PPW conjecture for domains in hemispheres of 
$\SS^n$: 

\begin{thm} \label{thm:I4}
Let $\Omega$ be contained in a hemisphere of $\SS^n$.  Then
$$
\lambda_2 / \lambda_1 (\Omega) \leq \lambda_2/ \lambda_1(\Omega^{\star})
$$
with equality if and only if $\Omega$ is itself a geodesic ball in $\SS^n$.
\end{thm} 

This theorem follows from Theorems \ref{thm:I1} and \ref{thm:I3}
above. If $B_{\lambda_1}$ is a ball having the same $\lambda_1$ as
$\Omega$ we have 
$$
\lambda_2 / \lambda_1 (\Omega) \leq \lambda_2 / \lambda_1 (B_{\lambda_1}) =
\lambda_2(\theta_1)/ \lambda_1(\theta_1),
$$ 
where $\theta_1$ is the geodesic radius of $B_{\lambda_1}$.  Theorem 
\ref{thm:I4} now follows from Theorem \ref{thm:I3} and Sperner's optimal 
Faber--Krahn-type result \cite{Sp73} for domains in $\SS^n$, which implies 
that $\theta_1(\Omega^{\star}) \ge \theta_1 (B_{\lambda_1})$ (via domain 
monotonicity). 

The results given in Theorems \ref{thm:I1}--\ref{thm:I4} were announced
earlier in \cite{AB94b} and \cite{AB97} (note, however, that Theorem 
\ref{thm:I2} is not mentioned explicitly in \cite{AB94b}, nor is Theorem 
\ref{thm:I3} stated formally there).  Moreover, \cite{AB97} contains 
alternative proofs of some of the results proved here.  These may be of 
independent interest, since they may point the way to generalizing the 
results presented here to other settings.  We note, though, that the 
proofs contained herein are those by which we first proved the results 
stated above.  
In fact, there is a close parallel between the proofs given here and 
those of our ``second'' proof of the Euclidean PPW conjecture given in 
\cite{AB92b} and of our $\SS^n$--version of the Szeg\H{o}--Weinberger 
inequality given in \cite{AB95}.  One of our objectives in \cite{AB94b} 
was to bring these similarities to the fore.

Beyond this, and as alluded to in passing in the previous paragraph, one 
might ask to what extent our results above are optimal, and, further, whether 
or not they are amenable to generalization (for example, to other spaces, 
the most obvious being $\HH^n$).  In particular, one might ask if it might 
be possible to remove our restriction to domains {\it contained in a 
hemisphere of $\SS^n$}.  It will be apparent to anyone who studies our
proofs that the success of our approach is very much dependent 
on this condition.  Certainly it should take little convincing to see that
geometrically things become quite different ``beyond the hemisphere'':  
for example, up to the hemisphere both the volume and surface area of a 
geodesic ball are increasing functions, but beyond the hemisphere the 
volume continues to grow while the surface area actually {\it shrinks}.  
On the other hand, we do not have any counterexamples to the conjecture 
that our theorems above (aside from Theorem \ref{thm:I2}, which is already 
established for all geodesic balls, and not just those contained in a 
hemisphere) continue to hold in the absence of the hemisphere condition.  
In fact, we have been able to establish that Theorem \ref{thm:I3} does 
hold for all $0<\theta_1<\pi$ for $n=2$ and $3$.  We might also note that 
it is common to encounter some sort of impediment at the hemisphere when 
dealing with the eigenvalues (specifically, Dirichlet or Neumann) of the 
Laplacian for domains in $\SS^n$.  For example, one might note that while 
for domains in Euclidean space and for domains strictly contained in a 
hemisphere of $\SS^n$ the first nonzero Neumann eigenvalue is always less 
than the first Dirichlet eigenvalue, the situation is reversed for 
geodesic balls which are larger than a hemisphere (equality obtains at 
the hemisphere).  For more discussion along these lines, the reader might 
consult \cite{AL}, \cite{AB94b}, \cite{AB95}, \cite{Ba80}, and references 
therein.

As for extensions to bounded domains contained in $\HH^n$, the
situation is still unsettled, but we make the following comments and 
observations.  The analog of Theorem \ref{thm:I4} (the ``naive'' PPW 
conjecture for $\HH^n$) cannot possibly hold, since one can imagine 
having a small disk (ball) with very narrow tentacles extending from it 
in such a way that the first two Dirichlet eigenvalues of this domain are 
very nearly those of the disk (ball), while the volume of the domain is 
as large as one wants.  This means that $\Omega^\star$ can be an
arbitrarily large ball, and the rub now comes from the fact that both
eigenvalues of the large ball can be made arbitrarily close to the bottom of 
the spectrum of the Laplacian on all of $\HH^n$, which is a positive value
(specifically $(n-1)^2/4$), by making the ball sufficiently large.  This
implies that $\lambda_2/\lambda_1 (\Omega^\star)$ can be made as close 
as we want to $1$, while $\lambda_2/\lambda_1 (\Omega)$ will stay close 
to its value for a small disk (ball).  Since this latter value is nearly 
the Euclidean value (any value larger than $1$ will do as well), 
we are faced with a contradiction.  In fact, it is likely that
$\lambda_2/\lambda_1$ for a geodesic ball in $\HH^n$ is a {\it decreasing}
function of the radius, i.e., that the counterpart to Theorem \ref{thm:I3} 
goes the other way (this certainly appears to be the case for $\HH^2$, 
based on numerical studies we have done).  This, of course, would be an 
interesting fact to prove in its own right (and perhaps especially for 
its geometric implications), even if it doesn't lead into a proof of the 
``PPW conjecture for $\HH^n$''.  It would also be, in a certain sense, 
the ``natural'' result, since $\lambda_2/\lambda_1$ is increasing for 
geodesic balls contained in a hemisphere of $\SS^n$ (and quite possibly 
for all geodesic balls in $\SS^n$; we note in this connection that
$\lambda_2/\lambda_1$ goes to infinity as the full sphere is approached 
since in that limit $\lambda_1 \rightarrow 0^+$ while $\lambda_2 
\rightarrow n^+$), while it is constant for all balls in $\RR^n$.  
As for an $\HH^n$-analog of Theorem \ref{thm:I1} (the ``sophisticated'' 
PPW conjecture for $\HH^n$), this may well be true, but as yet it is not 
proved.  One might also speculate that for bounded domains in $\HH^n$ 
$\lambda_2/\lambda_1$ is always less than the Euclidean bound (i.e., the 
value of $\lambda_2/\lambda_1$ for a Euclidean ball).  This conjecture is 
certainly supported by the behavior of $\lambda_2/\lambda_1$ for geodesic 
balls in $\HH^n$, and in general by the fact that ``large'' domains can be 
expected to have $\lambda_2/\lambda_1$ near $1$.  It is also supported by a 
result of Harrell and Michel \cite{HaMi94}, which gives a finite upper 
bound to $\lambda_2/\lambda_1$ for bounded domains contained in $\HH^2$.  
While their bound is almost certainly not optimal, it is at least in the 
right ballpark:  it is $17$, while the Euclidean value in two dimensions is 
approximately $2.5387$.

\bigskip

\section{The ``center of mass'' result for domains in spheres} 

\bigskip 

    In this section we present a general Center of Mass Theorem which 
follows from general topological arguments.  This Center of Mass
Theorem guarantees the orthogonality of certain functions which is
needed in later sections.  In a previous paper \cite {AB95} we 
gave such a theorem for domains contained in a hemisphere of $\SS^n$.
The proof of that theorem was fairly involved since we had to
identify a hemisphere of $\SS^n$ with the ball $B^n$ and use the
Brouwer fixed point theorem.  We also required a limiting argument
since there were certain problems which arose for  
domains extending to the equator that we could not handle directly. 
These arguments would suffice to yield a version of our Theorem
\ref{thm:CM1} below for domains contained in a hemisphere (see
Theorem 2.1 of \cite{AB95} and Remark 4 following it to see how to
make this extension).  Such a result would be enough to allow us to
prove only slightly weaker versions of the theorems found in the
remainder of this paper.  We have chosen to present the more general
result here since its proof is both simpler and more natural than our
former argument.  In the new proof attention is confined to mappings
from $\SS^n$ to $\SS^n$ and to certain natural geometrical
conditions.  A modest knowledge of degree theory is needed to conclude
the proof.  

\begin{thm} \label{thm:CM1} {\bf (Center of Mass Theorem)} 
Let $\Omega$ be a domain in $\SS^n$ and let $\tilde G$ be a
continuous function on  $[0,\pi]$ which is positive on $(0,\pi)$ and
symmetric about $\pi/2$.  Also, let $d \mu$ be any positive measure on 
$\Omega$.  Then there is a choice of Cartesian coordinates $x_1', x_2', 
\dots, x_{n+1}'$ for $\RR^{n+1}$ with the origin at the center of 
$\SS^n$ such that
\begin{equation}
\int_{\Omega} x_i' \, \tilde G(\theta) \, d\mu=0 \qquad 
\mbox{ for $i=1,2, \dots ,n$,}
\label{eq:CM1}
\end{equation}
where $\theta$ represents the angle from the positive $x_{n+1}'$--axis. 
\end{thm} 

\begin{rems}
In our later applications $d\mu$ will be $u_1^2 \, d \sigma$ where
$u_1$ is the first eigenfunction of $-\Delta$ on $\Omega$ with
Dirichlet boundary conditions imposed on $\partial \Omega$ and
$d\sigma$ is the standard volume element for $\SS^n$.  Similarly, in
our applications the function $\tilde G$ will be related to the
function $g(\theta)$ defined in equation (\ref{eq:4.3}) below 
(extended appropriately to $[0,\pi]$) by $\tilde
G(\theta)=g(\theta)/\sin \theta$.  Note that Theorem \ref{thm:CM1}
makes no statement about the integral in (\ref{eq:CM1}) for the case
$i=n+1$. 
\end{rems} 

\begin{proof}
We begin by considering the vector function $\vec v: \SS^n \to \RR^{n+1}$
defined for $\vec y \in \SS^n$ by
\begin{equation}
\vec v (\vec y) = \int_{\Omega} \vec x \, \tilde G(\theta) \,d\mu  
\label{eq:CM2}
\end{equation}
where $\theta$ represents the angle between $\vec y$ and the integration 
variable $\vec x \in \Omega$ and $\vec x=(x_1,x_2, \dots, x_{n+1})$
where $x_1, x_2, \dots, x_{n+1}$ represent some (initial) set of
Cartesian coordinates.  The vector $\vec v(\vec y)$ simply gives the
center of mass in $\RR^{n+1}$ of the hypersurface distribution on $\Omega 
\subset \SS^n$ with mass density given by $\tilde G(\theta) \, d\mu$.  
Note that its dependence on the point $\vec y$ is entirely through the 
function $\tilde G(\theta)$; indeed, if not for this function we would 
have only a single center of mass vector $\vec v$.  

First we argue that what we should look for are points $\vec y_0 \in 
\SS^n$ such that
\begin{equation}
\vec v (\vec y_0)= \alpha \vec y_0, 
\label{eq:CM3}
\end{equation}
i.e., such that $\vec y_0$ and $\vec v (\vec y_0)$ are linearly
dependent.  To see  that finding such a point will suffice to prove
our theorem, suppose $\vec y_0$ is a point where (\ref{eq:CM3})
holds and let $R$ be an $(n+1)$ by $(n+1)$ rotation matrix with 
$\vec y_0$ as its last row. Defining new Cartesian coordinates 
$x_1', x_2', \dots, x_{n+1}'$ via
\begin{equation}
x_i'= \sum_{j=1}^{n+1} R_{ij} x_j  
\label{eq:CM4}
\end{equation}
we have (with $\theta$ measured from $\vec y_0$)
\begin{eqnarray}
\int_{\Omega} x_i' \tilde G(\theta) \,  d\mu = \sum_{j=1}^{n+1} R_{ij} 
\int_{\Omega} x_j \tilde G(\theta) \, d\mu &=& {\left[R \vec v(\vec y_0) 
\right]}_{i}=\nonumber\\
&=&{\left[\alpha R\vec y_0 \right]}_{i}= {\left[ \alpha \hat e_{n+1} 
\right]}_{i}=0  \qquad \mbox{for $i=1,\dots,n$},  
\label{eq:CM5}
\end{eqnarray}
since $R$ is an orthogonal matrix and hence its first $n$ rows are
orthogonal to its last row.  Thus, the conclusion to the theorem will
follow from finding a  solution to (\ref{eq:CM3}) and from here on we
concentrate on finding such a solution.   

Now if $\vec v(\vec y)$ ever vanished for some $\vec y \in \SS^n$ the
conclusion  (\ref{eq:CM1}) would follow immediately with no need to
rotate coordinates.  So we may as well assume that $\vec v$ never 
vanishes on $\SS^n$.  Under this assumption we can pass to 
consideration of the mapping $\vec w: \SS^n \to \SS^n$ defined by
\begin{equation}
\vec w (\vec y)=\frac{\vec v(\vec y)}{\vert \vec v(\vec y)\vert}.
\label{eq:CM6} 
\end{equation}
The dependence condition (\ref{eq:CM3}) then reduces to 
\begin{equation}
\vec w(\vec y_0) = \pm \vec y_0
\label{eq:CM7}
\end{equation}
so that we are seeking a fixed point or an ``anti--fixed point'' of 
$\vec w$. Under the assumptions of the theorem it suffices to seek 
only fixed points, i.e., solutions to
\begin{equation}
\vec w(\vec y_0) =  \vec y_0.   
\label{eq:CM8}
\end{equation}
This follows from the symmetry of $\tilde G(\theta)$ about $\theta= \pi /2$ 
which implies that $\vec v (-\vec y)= \vec v(\vec y)$ for all $\vec y
\in \SS^n$ and thus that $\vec w (-\vec y)= \vec w (\vec y)$ for all
$\vec y \in \SS^n$.  Hence if $\vec y_0$ is a solution to
(\ref{eq:CM7}), then either $\vec y_0$ or $-\vec y_0$ 
must be a solution to (\ref{eq:CM8}) and we can therefore concentrate
solely on finding  solutions to (\ref{eq:CM8}).  (Geometrically, too,
it is more natural to view a point $\vec y_0$ satisfying
(\ref{eq:CM8}) as a center of mass of $\Omega$ than the point $- \vec
y_0$ satisfying $\vec w(- \vec y_0) = \vec y_0$, even though there is
no difference  between the two as far as fulfilling the conditions of
the theorem goes.)  

Finally we are at a point where we can use degree theory to conclude
that $\vec  w$ must have a fixed point, i.e., a solution $\vec y_0$ to
(\ref{eq:CM8}). The definition  
of $\vec v$ shows that it is continuous and our assumption that $\vec v$ 
does not vanish on $\SS^n$ guarantees the continuity of $\vec w$ as 
defined by (\ref{eq:CM6}).  Therefore, by standard theory (see, for 
example, \cite{Ar83}, p.~195; \cite{Du66}, Chapter 16, Section 1; 
\cite{HY61}, p.~263; \cite{Mu84}, p.~116; \cite{Sm71}, Chapter 16;
\cite{Wh63}, p.~807) $\vec w$ has a degree as a map from $\SS^n$ to $\SS^n$.
The fact that $\vec w(-\vec y)=\vec w(\vec y)$ for all $\vec y \in
\SS^n$ implies that $\deg (\vec w)$ is even (simply observe that for
an image point $\vec z \in \SS^n$ the points in the preimage come in
pairs  $\pm \vec y$).  Now if $\vec w$ has no fixed points, i.e., if
$\vec w (\vec y)  \neq \vec y$ for all $ \vec y \in \SS^n$, then
$\vec w$ would be homotopic to  the antipodal map $\vec A$ defined by
$\vec A(\vec y)=-\vec y$ via the homotopy 
\begin{equation}
h(t,\vec y)= \frac{-t\vec y + (1-t) \vec w(\vec y)}{\vert -t\vec y + (1-t) 
\vec w(\vec y) \vert} \qquad \mbox{for $0\le t \le 1$.}   
\label{eq:CM9}
\end{equation}
But the degree of a mapping is a homotopy invariant and $\deg (\vec A)= 
\pm 1 \neq \deg (\vec w)$, a contradiction.  Hence $\vec w$ must have some
fixed point  $\vec y_0 \in \SS^n$ and the proof is complete.
\end{proof}  

\begin{rems}  

(1) If $\Omega$ is contained in a hemisphere of $\SS^n$ then another way 
to complete the proof of the theorem is to set things up initially in a 
Cartesian frame such that $\Omega$ lies in the northern hemisphere of 
$\SS^n$ and observe that for any $\vec y \in \SS^n$
\begin{equation}
v_{n+1} (\vec y) = \int_{\Omega} x_{n+1} \tilde G(\theta) \, d\mu >0. 
\label{eq:CM10}
\end{equation}
Hence $\vec v$ never vanishes and we may regard $\vec w =\vec v/\vert 
\vec v \vert$ as a mapping from the closed northern hemisphere into itself.
Since this  space is homeomorphic to the ball $B^n$ we can apply the
Brouwer fixed point theorem to conclude that $\vec w$ has a fixed point 
$\vec y_0$ in the northern hemisphere (see, for example, \cite{Munk}, 
Section 8-10; \cite{Sm71}, p.~406; \cite{Sp66}, pp.~151, 194).  This proof 
is similar to, but simpler than, the proof of the restricted version of 
Theorem \ref{thm:CM1} that we gave in \cite{AB95}.  

(2) The alternative proof just given does not use the symmetry of $\tilde 
G(\theta)$ about $\theta=\pi/2$ stated in the theorem.  Only positivity of 
$\tilde G(\theta)$ on $(0,\pi)$ is used.  There are certainly situations, 
in particular for domains which are in some sense larger than a hemisphere, 
where one might not want to require that $\tilde G$ be symmetric about 
$\pi/2$.  In such situations one possible route to a center of mass result 
is to show that the mapping $\vec w$ misses at least one point of $\SS^n$, 
following the spirit of the alternative proof given above.  Brouwer's fixed
point theorem can then be applied (to $\SS^n$ less a sufficiently
small neighborhood of a point that $\vec w$ misses) to yield a fixed
point.  Or in the language of degree theory, the case where a map
misses a point is the simplest case of a map which has degree $0$.  
This can be seen directly or by observing that such a map $\vec w$
is homotopic  to a constant map, i.e., contractible to a point, or
{\it inessential} (see \cite{HY61}, p.~154; \cite{Munk}, p.~357; or 
\cite{Sp66}, p.~23) via the homotopy 
\begin{equation}
g(t,\vec y)= \frac{-t\vec z_0 + (1-t) \vec w(\vec y)}{\vert -t\vec z_0 + 
(1-t) \vec w(\vec y) \vert} \qquad \mbox{for $0\le t \le 1$}   
\label{eq:CM11}
\end{equation}
if $\vec z_0$ is a point that $\vec w$ misses.   

(3) If $\Omega$ has a center of symmetry and $d \mu$ is either $u_1^2 \, 
d \sigma$ (as occurs in Dirichlet problems for $-\Delta$ on
$\Omega$) or $d \sigma$ (as occurs in Neumann problems for $-\Delta$
on $\Omega$), then this point will certainly serve as a center of mass 
in the sense of Theorem \ref{thm:CM1}.  Here $d\sigma$ represents the 
standard volume element for $\SS^n$ and $u_1$ represents the  first
Dirichlet eigenfunction of $-\Delta$ on $\Omega$.  One has only to note 
that $u_1$, being unique up to a constant factor and of one sign, 
must share the symmetries of $\Omega$.  Also, if $d\mu=d\sigma$ and 
$\Omega$ is an arbitrary domain such that its complement has a center 
of mass in the sense of Theorem \ref{thm:CM1} (as could be concluded, 
for example, via any of the conditions discussed so far, or by other 
means) then $\Omega$ shares this center of mass since it is clear that 
\begin{equation}
\int_{\SS^n} x_i' \tilde G(\theta) \, d \sigma=0 \qquad \mbox{for $i=1,2 
\dots,n$}
\end{equation}
if $x_1', \dots ,x_{n+1}'$ represent Cartesian coordinates and $\theta$ 
is measured from the positive $x_{n+1}'$--axis.  In particular, 
$\Omega$ certainly has a center of mass in this sense if the complement 
of $\Omega$ is contained in a hemisphere.  All the observations made in
this remark apply whether or not $\tilde G(\theta)$ is even with respect 
to $\theta=\pi/2$.
For another result that holds in the absence of symmetry of $\tilde G$ about 
$\pi/2$ see Theorem \ref{thm:CM2} below.  

(4) With Theorem \ref{thm:CM1} in hand we can obtain a modest improvement 
of our main theorem in \cite{AB95} (see Theorem 5.1).  In particular, in 
Remark 2 following the proof of Theorem 5.1 in \cite{AB95} we now have no 
need to invoke condition (2.4).  We thus obtain the result $\mu_1(\Omega)
\le \mu_1(\Omega^{\star})$ comparing the first nonzero Neumann eigenvalue
of the domain $\Omega$ with that of the spherical cap $\Omega^{\star}$  
having the same volume for any domain $\Omega$ such that $\Omega \cap 
(-\Omega)=\emptyset$ (equivalently, $-\Omega \subset \SS^n\setminus \Omega$). 
This inequality is an equality if and only if $\Omega$ is itself a 
geodesic ball.  
More generally, the same result holds if, when the 
north pole is a center of mass in the sense of Theorem \ref{thm:CM1}, 
$\Omega$ has the property that for each $k \in  [0,1]$
\begin{equation}
\vert \{ \vec y \in \Omega \bigm| y_{n+1} \le -k \} \vert \le 
\vert \{ \vec y \in \SS^n \setminus \Omega \bigm| y_{n+1} \ge k \} \vert   
\label{eq:CM12}
\end{equation}
(cf.\ Remark 2 following the proof of Theorem 5.1 in \cite {AB95}).  Here 
$\vert X \vert$ denotes the measure of $X$ with respect to the canonical 
measure (standard volume element) on $\SS^n$ where $X$ is any measurable 
subset of $\SS^n$.  We refer to condition (\ref{eq:CM12}) as the ``excess 
less than or equal to deficit property''.  For some further comments 
relating to this property, see our remarks at the end of Section 6.  These 
give the most general conditions known at this time.   
\end{rems} 

Finally, for possible future use (see also Remarks 2 and 3 above) we 
state the following variant of our Center of Mass Theorem which holds in 
even dimension (i.e., for $\Omega \subset \SS^n$ with $n$ even) in 
the absence of symmetry of $\tilde G$ about $\pi/2$: 

\begin{thm}  \label{thm:CM2}
Let $\Omega$ be a domain in $\SS^n$ for $n$ even and let
$\tilde G$  be continuous on $[0,\pi]$ and positive on $(0,\pi)$.
Then $\Omega$ has a center of mass in the sense of Theorem \ref{thm:CM1}. 
That is, there exists a choice of Cartesian coordinates such that
(\ref{eq:CM1}) holds for $i=1,2, \dots, n$.  
\end{thm} 

\begin{proof} 
Defining $\vec v$ as above (equation (\ref{eq:CM2})) either $\vec v$
vanishes somewhere and we are done or we can pass to $\vec w=\vec
v/\vert \vec v \vert$.  Continuing with the latter case, if $\vec w$
has neither a fixed point nor an anti--fixed  point (i.e., there are
no solutions $\vec y \in \SS^n$ to $\vec w(\vec y)=\pm  \vec y$) then
as above we can show that $\vec w$ is homotopic to the antipodal map
and also to the identity map.  But the antipodal map has degree 
$(-1)^{n+1}$  (see, for example, \cite{Ar83}, p.~197, Theorem 9.2; 
\cite{Du66}, p.~339, Exercise 4; \cite{Mu84}, p.~118, Theorem 21.3;
\cite{Sm71}, p.~403, Theorem 4.3; \cite{Wh63}, p.~809) and for $n$ 
even $(-1)^{n+1}=-1 \neq 1 =$ degree of the identity map.  This is 
a contradiction since degree is a homotopy invariant (see \cite{Ar83}, 
p.~195; \cite{Du66}, p.~339; \cite{HY61}, p.~266; \cite{Mu84}, p.~117;
\cite{Sm71}, p.~401; or \cite{Wh63}, p.~809), hence $\vec w$ must have 
either a fixed point or an anti--fixed  point, and the conclusion 
of the theorem follows.  
\end{proof} 

\bigskip

\section{Properties of the first two Dirichlet eigenvalues of
geodesic balls in $\SS^n$} 

\bigskip

    We consider the Dirichlet problem on a geodesic ball of radius
$\theta_1$ in $\SS^n$ (where $\theta_1 \in (0,\pi)$) which we view as 
a polar cap centered at the
north pole (i.e., the point $\hat e_{n+1} \in \SS^n$ is taken as the
center of our geodesic ball). By using the $O(n)$ symmetry of the
polar cap (this $O(n)$ is the subgroup of $O(n+1)$ which leaves the
point $\hat e_{n+1}$ fixed), one can separate variables in the usual
way obtaining the family of ordinary differential equations in the 
``radial'' variable $\theta$
\begin{equation}
-y''-(n-1) \cot \theta \, y'+m(m+n-2)\csc^2 \theta \, y=\lambda y \qquad
\mbox{on $(0,\theta_1)$}
\label{eq:3.1}
\end{equation}
for $m=0,1,2,\dots$. The boundary conditions to be applied for
(\ref{eq:3.1}) are $y(0)$ finite and $y(\theta_1)=0$.  In particular,
we shall be concerned with the lowest eigenvalues of the $m=0$ and
$m=1$ cases of this equation, but first we develop some general
properties of the solutions to these equations and some of their
interrelationships. 
  
    We begin by considering $\lambda$ as a positive parameter (all the 
eigenvalues that we consider here are easily seen to be positive by 
consideration of the Rayleigh quotients that characterize them) and 
defining $u_m(\theta;\lambda)$  for $m=0,1,2, \dots$ as that solution to
(\ref{eq:3.1}) which has the behavior
\begin{equation}
u_m(\theta;\lambda)=c_m \theta^m + O(\theta^{m+2})
\label{eq:3.2}
\end{equation}
where the constants $c_m$ will be specified below.  This behavior is
consistent with  equation (\ref{eq:3.1}) as can be seen from
Frobenius theory ($\theta=0$ is a regular singular point of
(\ref{eq:3.1})).  In particular, the eigenfunctions to 
(\ref{eq:3.1}) will all be found among the $u_m$'s defined in 
(\ref{eq:3.2}) assuming $c_m\neq 0$ since finiteness at $\theta=0$
forces this behavior (up to multiplicative factors).  Moreover, it is
not difficult to verify that if $u_m$ solves 
(\ref{eq:3.1}) then
$$
u_m'-m \cot \theta \,  u_m
$$
satisfies (\ref{eq:3.1}) for $m$ replaced by $m+1$ and also
$$
u_m'+(m+n-2)\cot \theta \, u_m
$$
satisfies (\ref{eq:3.1}) for $m$ replaced by $m-1$.  From these facts
and Frobenius theory again it follows that
\begin{equation}
u_{m+1} = -u_m'+m \cot \theta \, u_m
\label{eq:3.3}
\end{equation}
and
\begin{equation}
\left[\lambda -(m-1)(m+n-2) \right] u_{m-1}=u_m'+(m+n-2)\cot \theta
\, u_m
\label{eq:3.4}
\end{equation}
if we agree to set $c_0=1$ and define successive $c_m$'s via
\begin{equation}
c_{m+1}=\frac{\lambda-m(m+n-1)}{2m+n} c_m \qquad \mbox{for $m=0,1,2, \dots$.}
\label{eq:3.5}
\end{equation}
These constitute the raising and lowering relations for the functions
$u_m$.  Also, elimination of $u_m'$ between (\ref{eq:3.3}) and
(\ref{eq:3.4}) yields the pure recursion relation in $m$
\begin{equation}
u_{m+1} - (2m+n-2) \cot \theta \, u_m +
\left[\lambda-(m-1)(m+n-2)\right] u_{m-1}=0.
\label{eq:3.6}
\end{equation}
Since $(\sin \theta)^{m+n-2} u_m(\theta)$ is $0$ at $\theta=0$ one
can integrate (\ref{eq:3.4}) and replace $m$ by $m+1$ to obtain
\begin{equation}
(\sin \theta)^{m+n-1} u_{m+1}(\theta)=\left[\lambda-m(m+n-1)\right]
\int_0^{\theta} (\sin t)^{m+n-1} u_m(t) \, dt.
\label{eq:3.7}
\end{equation}
What we have developed so far could be considered the algebraic
properties of the functions $u_m$.  One should think of the $u_m$'s as
the analogs of associated Legendre functions in $n$ dimensions.  In
particular, when $n=2$ one finds that $u_m(\theta;\nu(\nu+1))
=(-1)^m P_{\nu}^m(\cos \theta)$ following \cite{AS64} (in general this 
should be correct up to a constant factor, typically a factor of 
$(-1)^m$, depending on the precise convention adopted; cf.\ \cite{Ri81}, 
for example, where the convention differs by $(-1)^m$); here we have 
replaced the eigenvalue parameter $\lambda$ by $\nu(\nu+1)$, as is 
traditional in dealing with Legendre functions.  In this case 
(\ref{eq:3.1}) reduces to the associated Legendre equation and 
(\ref{eq:3.3}), (\ref{eq:3.4}), and (\ref{eq:3.6}) all reduce to standard 
relations between the associated Legendre functions.  Almost certainly 
this generalization to $n$ dimensions of associated Legendre functions 
and their basic relations has been developed before, though we do not 
know of a reference where the details needed here are developed explicitly.  
(Cf.\ also \cite{FrHa76} which deals with the $m=0$ case of (\ref{eq:3.1}) 
in $n$ dimensions but in Schr\"{o}dinger normal form.)  In any event, it 
is a relatively simple matter to pass to the $n$--dimensional case once 
the situation in two dimensions is known.  In fact, our functions $u_m$ 
can be expressed in terms of associated Legendre functions no matter what 
the dimension but this connection is not particularly useful in this 
context so we do not elaborate upon it here (but see the equations in 
Section~6 following equation (\ref{eq:6.4}), or Remark $1$ below).
We note that our convention on the constants $c_m$ allows them to vanish 
from a certain value of $m$ on for specific values of the parameter 
$\lambda$ (cf.\ (\ref{eq:3.5})).  (Recall, for example, that when dealing 
with the full sphere $\SS^2$ one needs $P_{\ell}^m$ only for $\ell=0,1,
\dots$ and $m=0,1,\dots, \ell$.)  This is not a problem for our purposes 
here since we will be most interested in passing from $u_m$ to $u_{m+1}$ 
via either (\ref{eq:3.3}) or (\ref{eq:3.7}).  

We come back now to the eigenvalue problem for (\ref{eq:3.1}).  If one
defines the left-hand side of equation (\ref{eq:3.1}) as the operator
$h_m$ applied to $y$ (with boundary conditions incorporated in the
definition of $h_m$), then it is easily seen that $h_{m'}>h_m$ in the
sense of quadratic forms if $m'>m$. Thus $\lambda_1$ of $-\Delta$ for
the polar cap is $\lambda_1(h_0)$ while $\lambda_2$ of $-\Delta$ for
the polar cap must be either $\lambda_1(h_1)$ or $\lambda_2(h_0)$.  We
now show that the former is the case. 

\begin{lemma} \label{lem:3.1}
The first eigenvalue of the Dirichlet Laplacian on a
polar cap is the first eigenvalue of (\ref{eq:3.1}) with $m=0$ while
the second eigenvalue of the Dirichlet Laplacian is the first
eigenvalue of (\ref{eq:3.1}) with $m=1$.  The second eigenvalue of the
cap occurs with multiplicity $n$.  These results hold for all polar caps, 
i.e., for all $\theta_1 \in (0, \pi)$.
\end{lemma}

\begin{proof}
In the notation of the preceding paragraph we must show that
$\lambda_2(h_0) >\lambda_1(h_1)$.  The argument proceeds via a simple
use of Rolle's theorem as applied to (\ref{eq:3.3}) and (\ref{eq:3.4})
rewritten in the forms
\begin{equation}
\left[(\sin \theta)^{-m} u_m \right]'= -(\sin \theta)^{-m} u_{m+1}
\label{eq:3.8}
\end{equation}
and
\begin{equation}
\left[(\sin \theta)^{m+n-2} u_m \right]'= \left[\lambda-(m-1)(m+n-2)
\right] (\sin \theta)^{m+n-2} u_{m-1}.
\label{eq:3.9}
\end{equation}
In particular, with $m=0$ in the first of these we have
\begin{equation}
u_0'=-u_1
\label{eq:3.10}
\end{equation}
and with $m=1$ in the second we have
\begin{equation}
\left[(\sin \theta)^{n-1} u_1 \right]'= \lambda(\sin \theta)^{n-1} u_{0}.
\label{eq:3.11}
\end{equation}
By Rolle's theorem, between any two zeros of $u_0$ there is a zero of
$u_0'$ and hence of $u_1$, since (\ref{eq:3.10}) holds.  Similarly,
between two zeros of $(\sin \theta)^{n-1} u_1$ there is a zero of its
derivative and hence, by (\ref{eq:3.11}), of $u_0$.  Thus for fixed 
$\lambda>0$ the zeros of $u_0$ and $(\sin \theta)^{n-1} \, u_1$ on 
$[0,\pi)$ interlace.  

Now consider $u_0$ and $u_1$ for $\lambda=\lambda_1(h_1)$. Since this
makes $\theta_1$ the first positive zero of $u_1$ it is clear by what
we have just proved that $u_0$ has exactly one zero in
$(0,\theta_1)$ and that $\theta_1$ is not a zero of $u_0$.  This then
implies, by the fact that the positive zeros of any $u_m$ are decreasing 
functions of the parameter $\lambda$ (see, for example, \cite{BR}, 
p.~315, or \cite{CoHi}, p.~454), that $\lambda_2(h_0)>\lambda_1(h_1)$. 

That the multiplicity of $\lambda_1(h_1)$ as an eigenvalue of
$-\Delta$ on the polar cap (= geodesic ball) with Dirichlet boundary
conditions is $n$ follows from the details of separation of
variables. It can be shown that the corresponding eigenfunctions can
be taken as $(x_i /\sin \theta) y (\theta)$ (restricted to $\SS^n$)
where $i=1,2, \dots, n$, $x_{n+1} = \cos \theta$, and $y(\theta)$ is
the eigenfunction of (\ref{eq:3.1}) for the eigenvalue
$\lambda_1(h_1)$.  These functions form an orthogonal basis for the
eigenspace of $-\Delta$ corresponding to the eigenvalue
$\lambda_1(h_1)$, showing that its multiplicity is $n$. This
completes our proof. 
\end{proof} 

\begin{rems}

(1) The argument used in our proof above can be viewed as a way of 
translating order properties of the zeros of the $u_m$'s into order 
properties of the corresponding eigenvalues $\lambda_i(h_m)$.  In fact, 
our approach extends easily to an interlacing result for the zeros 
(and hence for the associated eigenvalues) of $u_m$ and $u_{m+1}$ for 
arbitrary $m$.  Further ordering properties of the Dirichlet eigenvalues 
of spherical caps in $\SS^n$, at least for even $n$ (and surely the case 
of odd $n$ could be handled similarly), may be inferred from the papers 
of Baginski \cite{Bag1}, \cite{Bag2}.  We note, though, that all of 
Baginski's results are presented in terms of the zeros of the associated 
Legendre functions $P_{\nu}^m$ in the variable $\nu$ (for $m$ an integer).  
To make the connection to the present setting, one makes use of the 
formulas in Section 6 following equation (\ref{eq:6.4}) (or equations 
(3.17) and (3.18) in \cite{AB94b}, with the correction that the upper 
index in both associated Legendre functions should be negated), which 
relate our functions $u_m$ as defined above to associated Legendre 
functions (up to constant factors).  Note, in this connection, that 
$y_1(\theta)=u_0(\theta;\lambda_1)$ and $y_2(\theta)=u_1(\theta;\lambda_2)$ 
up to constant factors (as proved in Lemma \ref{lem:3.1} above; throughout 
this section we take these factors to be $1$).  In general, one has 
$u_m(\theta;\lambda)$ proportional to $(\sin \theta)^{1-n/2} \, 
P_{\nu}^{-(n/2-1+m)}(\cos \theta)$, where $\lambda$ and $\nu$ are 
related by $\lambda=(\nu-n/2+1)(\nu+n/2)$.  We thank the referee for 
calling Baginski's papers to our attention.     

(2) Another proof of Lemma \ref{lem:3.1} follows by mimicking our
proof of Lemma 3.1 of \cite{AB95}.  With $\tau=\lambda_2(h_0)$ and $v$ 
as the associated eigenfunction we can assume that for some $a\in
(0,\theta_1)$ $v>0$ on $(0,a)$ and $v<0$ on $(a,\theta_1)$, and hence
that $v(a)=0$, $v'(a)<0$.  Also we take $\lambda=\lambda_1(h_1)$ and
set $g=u_1(\theta;\lambda)$ and $h=u_0(\theta;\lambda)$. Since
$c_1=\lambda/n>0$ it is clear that $g>0$ on $(0,\theta_1)$ and that
$g(0)=0=g(\theta_1)$. It also follows that  $g=-h'$ and $h$ satisfies
$$
-h''-(n-1)\cot \theta \, h'=\lambda h
$$
while $v$ satisfies
$$
-v''-(n-1)\cot \theta \, v'=\tau v.
$$
From the last two equations we obtain
$$
\left[(\sin \theta)^{n-1}(vh'-v'h)\right]'=-(\lambda-\tau)(\sin
\theta)^{n-1} v h 
$$
and integration from $a$ to $\theta_1$ produces
\begin{equation}
(\sin \theta_1)^{n-1} v'(\theta_1)h(\theta_1) -(\sin a)^{n-1}
v'(a)h(a)=(\lambda-\tau)\int_a^{\theta_1} v h \sin^{n-1} \theta \, d
\theta. 
\label{eq:3.12}
\end{equation}
We now argue by contradiction, so assume $\tau \le \lambda$. 
Since $v(a)=0$, $v=u_0(\theta;\tau)$, $h=u_0(\theta;\lambda)$ and
$\tau \le \lambda$, by the fact that the positive zeros of
$u_0(\theta;\lambda)$ are decreasing with increasing $\lambda$ it 
follows that the first positive zero of $h$ is less than or equal to 
$a$ and, since $h'=-g<0$ on $(0,\theta_1)$, it therefore follows that 
$h<0$ on $(a,\theta_1]$.  But now (\ref{eq:3.12}) gives a contradiction, 
since its right-hand side is greater than or equal to $0$ while its
left-hand side is negative (note that $v'(a)<0$ and $v'(\theta_1)>0$
since $a$ and $\theta_1$ must be successive zeros of $v$).   

(3) Yet another proof of Lemma \ref{lem:3.1} would be via the 
level--ordering results of Baumgartner, Grosse, and Martin \cite{BGM84}, 
\cite{BGM85}.  Specifically, see our papers \cite{AB88}, \cite{AB88b} 
where proofs for a ball in $\RR^n$ occur and also the papers \cite{Bau86}, 
\cite{Bau92} of Baumgartner, which give extensions to cases arising from 
separation of variables in spherical coordinates in spaces of constant 
curvature.
\end{rems} 

\bigskip

For future reference we note the following lemma. 

\begin{lemma} \label{lem:3.2}
If $0 < \theta_1 < \pi$, the first eigenfunction of (\ref{eq:3.1})
with $m=0$, i.e., $y_1(\theta) \equiv u_0(\theta;\lambda_1)$, is
strictly decreasing on $[0,\theta_1]$ ($y_1 > 0$ on $[0,\theta_1)$ 
is our convention for $y_1$ here and throughout this paper; this 
follows from our choice $c_0=1$). 
\end{lemma}

\begin{proof} 
$u_0(\theta;\lambda_1)$ satisfies
$$
- (\sin^{n-1} \theta \, u_0')'=\lambda_1(\theta_1) \,
\sin^{n-1} \theta \, u_0 >0
$$
in $[0,\theta_1)$ (since $\lambda_1(\Omega) > 0$ follows from the 
variational characterization of the eigenvalues of $-\Delta$ via the 
Rayleigh quotient $\int_\Omega |\nabla \varphi|^2 /\int_\Omega \varphi^2$), 
which implies that $\sin^{n-1} \theta \, u_0'$ is decreasing in 
$[0,\theta_1)$.  Hence 
$\sin^{n-1} \theta \, u_0'<(\sin^{n-1} \theta \, u_0')\bigm|_{\theta=0}=0$,
which proves the lemma. 
\end{proof} 
   
\bigskip

Having identified $\lambda_1$ and $\lambda_2$ for $-\Delta$ on a
spherical cap with Dirichlet boundary conditions we are now in
position to investigate their behaviors and, in particular, that of
$\lambda_2/\lambda_1$.  Since our concern will be with how these
functions vary with $\theta_1$, the geodesic radius of the spherical
cap, we shall denote $\lambda_1$ and $\lambda_2$ by
$\lambda_1(\theta_1)$ and $\lambda_2(\theta_1)$ throughout the
remainder of this section. 
Associated with equation (\ref{eq:3.1}) is the one--dimensional
Schr\"odinger operator
\begin{equation}
H_{m}(\theta_1)=-\frac{d^2}{d\theta^2}+\frac{(2m+n-1)(2m+n-3)}
{4\sin^2\theta} -\frac{(n-1)^2}{4}
\label{eq:3.15}
\end{equation}
acting on $L^2((0,\theta_1),d\theta)$ with Dirichlet boundary 
conditions imposed at $0$ and $\theta_1$.  The operators
$H_m(\theta_1)$ form a family of self--adjoint operators.  It is clear 
by Lemma \ref{lem:3.1} that $\lambda_1(\theta_1)=\lambda_1(H_0(\theta_1))$ 
and $\lambda_2(\theta_1)=\lambda_1(H_1(\theta_1))$. 

We now analyze how $\lambda_1$ and $\lambda_2$ vary with $\theta_1$
by using perturbation theory \cite{HiSi96}, \cite{Ka76}, \cite{ReSi78}. 
To be successful at this we need to work on a fixed interval 
$(0,\theta_1)$ and we do this by observing that the eigenvalue problem
$H_m(c\theta_1)v=\lambda v$ on $(0,c\theta_1)$ can be rescaled to
\begin{equation}
\left[-\frac{1}{c^2}\frac{d^2}{dt^2}+\frac{(2m+n-1)(2m+n-3)}{4\sin^2
ct}-\frac{(n-1)^2}{4} \right]v=\lambda v \qquad 
\mbox{for $t\in (0,\theta_1)$}
\label{eq:3.17}
\end{equation}
which is equivalent to
\begin{equation}
\left[-\frac{d^2}{d\theta^2}+\frac{(2m+n-1)(2m+n-3)c^2}{4\sin^2
c\theta}-\frac{(n-1)^2c^2}{4} \right]v=c^2\lambda v \qquad \mbox{for
$\theta \in (0,\theta_1)$}.
\label{eq:3.18}
\end{equation}
As was done above for $H_m(\theta_1)$, we define an operator
$\tilde H_m(c)$ on $L^2(0,\theta_1)$ via the differential expression
appearing on the left-hand side of (\ref{eq:3.18}).  It is clear from
(\ref{eq:3.18}) that $\lambda_k(\tilde H_m(c))=c^2\lambda_k
(H_m(c\theta_1))$ and thus, in particular, that
\begin{equation}
\lambda_1(c\theta_1)=c^{-2}\lambda_1(\tilde H_0(c))
\label{eq:3.19}
\end{equation}
and
\begin{equation}
\lambda_2(c\theta_1)=c^{-2}\lambda_1(\tilde H_1(c)). 
\label{eq:3.20}
\end{equation} 

What we intend to do now is to determine the derivatives
$\lambda_1'(\theta_1)$ and $\lambda_2'(\theta_1)$ using perturbation
theory and the fact that
\begin{equation}
\lambda_j'(\theta_1)=\frac{1}{\theta_1}
\frac{d\lambda_j(c\theta_1)}{dc}\Bigm|_{c=1}. 
\label{eq:3.21}
\end{equation}
Since $\tilde H_m(c)$ is an analytic family in $c$ for $c$ near $1$
we can apply regular Rayleigh--Schr\"odinger 
perturbation theory \cite{HiSi96}, \cite{Ka76}, \cite{ReSi78}. 
In fact
\begin{eqnarray}
\tilde H_m(c) = H_m(\theta_1)&+&\frac{(2m+n-1)(2m+n-3)}{4}
\left(\frac{c^2}{\sin^2 c\theta}-\frac{1}{\sin^2 \theta} \right)
-\frac{(n-1)^2(c^2-1)}{4} 
\nonumber\\
&=& H_m(\theta_1)+V_m(\theta;c) \label{eq:3.22} 
\end{eqnarray}
and since $V_m(\theta;c)$ is analytic in $c$ for $c$ near
$1$ we are assured that the operators $\tilde H_m(c)$ form an
analytic family of type (A) for $c$ near 1 (see \cite{ReSi78}, p.~16 
for the definition of analytic family of type (A), or see
\cite{HiSi96}, p.~154; also see Chapter~7 of \cite{Ka76} for the
definitive account of analytic perturbation theory).  This allows us to
compute the derivatives of the eigenvalues of $\tilde H_m(c)$ using the 
first-order perturbation formula (cf.\ Kato \cite{Ka76}, p.~391, eq.~(3.18))
\begin{equation}
\frac{d\lambda_1(\tilde H_m(c))}{dc}\Bigm|_{c=1}=\frac{\int_0^{\theta_1}
\left[\frac{\partial V_m}{\partial c}(\theta;c)\bigm|_{c=1}
\right]v_m^2 \, d\theta}{\int_0^{\theta_1} v_m^2 \, d\theta}
\label{eq:3.23}
\end{equation}
where the functions $v_m(\theta)$ denote first eigenfunctions of
$H_m(\theta_1)=\tilde H_m(1)$. 

By (\ref{eq:3.21}) we have
\begin{eqnarray}
\frac{d}{d\theta_1} \left[
\frac{\lambda_2(\theta_1)}{\lambda_1(\theta_1)} \right] &=& 
\frac{1}{\theta_1} \left[ \frac{d}{dc} \left( 
\frac{\lambda_2(c \theta_1)}{\lambda_1(c \theta_1)} \right) \right]
\Bigm|_{c=1} 
\nonumber \\
&=&\frac{1}{\theta_1} \left[ \lambda_1 (\theta_1) \left( \frac{d}{dc}
\lambda_2(c \theta_1)\Bigm|_{c=1} \right) - 
\lambda_2 (\theta_1) \left( \frac{d}{dc}
\lambda_1 (c \theta_1) \Bigm|_{c=1} \right) \right]
\frac{1}{\lambda_1(\theta_1)^2}. 
\nonumber
\end{eqnarray}
Thus, showing that $\lambda_2(\theta_1)/\lambda_1(\theta_1)$
increases with increasing $\theta_1$ comes down to showing that
\begin{equation}
0<\frac{1}{\lambda_2(\theta_1)}
\left(\frac{d}{dc}\lambda_2(c \theta_1)\Bigm|_{c=1}\right)- 
\frac{1}{\lambda_1(\theta_1)}
\left(\frac{d}{dc}\lambda_1(c \theta_1)\Bigm|_{c=1}\right), 
\label{eq:3.24}
\end{equation}
which, by equations (\ref{eq:3.19}) and (\ref{eq:3.20}), reduces to
showing 
\begin{equation}
0<\frac{1}{\lambda_2(\theta_1)}
\left(\frac{d}{dc}\lambda_1(\tilde H_1(c))\Bigm|_{c=1}\right)- 
\frac{1}{\lambda_1(\theta_1)}
\left(\frac{d}{dc}\lambda_1(\tilde H_0(c))\Bigm|_{c=1}\right). 
\label{eq:3.25}
\end{equation} 

From (\ref{eq:3.22}) we obtain
\begin{equation}
\frac{\partial V_m}{\partial c} (\theta;c) \Bigm|_{c=1}
=\frac{1}{2}(2m+n-1)(2m+n-3)\csc^2 \theta (1-\theta \cot
\theta)-\frac{1}{2}(n-1)^2 
\nonumber
\end{equation}
so that, by (\ref{eq:3.23}),
\begin{equation}
\frac{d}{dc} \lambda_1(\tilde H_1(c)) \Bigm|_{c=1} = \frac{1}{2}(n-1)
\frac{\int_0^{\theta_1} \left[(n+1)\csc^2 \theta (1-\theta \cot
\theta) - (n-1) \right] v_1^2 \, d\theta}{\int_0^{\theta_1} v_1^2 \,
d \theta},
\label{eq:3.26}
\end{equation}
and
\begin{equation}
\frac{d}{dc} \lambda_1(\tilde H_0(c)) \Bigm|_{c=1} = \frac{1}{2}(n-1)
\frac{\int_0^{\theta_1} \left[(n-3)\csc^2 \theta (1-\theta \cot
\theta) - (n-1) \right] v_0^2 \, d\theta}{\int_0^{\theta_1} v_0^2 \,
d \theta}.
\label{eq:3.27}
\end{equation}
The functions $v_0$ and $v_1$ are related to $u_0(\theta;\lambda_1)$
and $u_1(\theta;\lambda_2)$ by $v_0=u_0 \sin^{(n-1)/2}\theta$ and 
$v_1=u_1 \sin^{(n-1)/2}\theta$ respectively. 
  
Introducing the functions $\ell(\theta)=\cot \theta - \theta \,
\csc^2 \theta$ and $m(\theta)= -\ell'(\theta)/2=\csc^2 \theta \,
(1-\theta\, \cot \theta)$ we can write (\ref{eq:3.26}) as
\begin{equation}
\frac{d}{dc} \lambda_1(\tilde H_1(c)) \Bigm|_{c=1} = \frac{1}{2}(n-1)
\frac{\int_0^{\theta_1} \left[(n+1) m(\theta) - (n-1) \right] v_1^2
\, d\theta}{\int_0^{\theta_1} v_1^2 \, d \theta},
\label{eq:3.28}
\end{equation}
and (\ref{eq:3.27}) as
\begin{equation}
\frac{d}{dc} \lambda_1(\tilde H_0(c)) \Bigm|_{c=1} = \frac{1}{2}(n-1)
\frac{\int_0^{\theta_1} \left[(n-3) m(\theta) - (n-1) \right] v_0^2
\, d\theta}{\int_0^{\theta_1} v_0^2 \, d \theta},
\label{eq:3.29}
\end{equation}
respectively.  Using the relations $m(\theta)=-\ell'(\theta)/2$, 
$\ell(\theta) \cot \theta =m(\theta)-1$, and $v_0^2= u_0^2 \sin^{n-1}
\theta$, we can rewrite the numerator
of the right-hand side of (\ref{eq:3.29}) as
\begin{eqnarray}
\int_0^{\theta_1} [(n-1)(m(\theta) - 1) -2 \, m(\theta)] \, v_0^2 
\, d \theta 
&=& \int_0^{\theta_1} [(n-1) \ell (\theta) \cot \theta +
\ell'(\theta)] \, u_0^2 \, \sin^{n-1} \theta \, d \theta
\nonumber\\
&=& \int_0^{\theta_1}[\ell(\theta) \, \sin^{n-1} \theta]' \, u_0^2 \, 
d \theta 
\nonumber\\
&=& \int_0^{\theta_1}[-2 \ell (\theta)\, u_0(\theta) \, u_0'(\theta)] 
\, \sin^{n-1} \theta \, d \theta, 
\label{eq:3.30}
\end{eqnarray}
where the last equality follows by integrating by parts (both boundary 
terms vanish).  Finally, from (\ref{eq:3.29}) and (\ref{eq:3.30}) 
we obtain
\begin{equation}
\frac{d}{dc} \lambda_1(\tilde H_0(c)) \Bigm|_{c=1} = (n-1)
\frac{\int_0^{\theta_1} \left[-\ell (\theta)\, u_0(\theta) \, 
u_0'(\theta) \right] \, \sin^{n-1}\theta \, d \theta}
{\int_0^{\theta_1} v_0^2 \, d \theta}.
\label{eq:3.31}
\end{equation} 

At this point we need the following properties of the functions 
$\ell(\theta)$ and $m(\theta)$.  

\begin{lemma} \label{lem:3.3}
The function
\begin{equation}
\ell(\theta) \equiv \cot \theta -\theta \csc^2 \theta
\label{eq:3.32}
\end{equation}
is negative, decreasing, and concave for $0<\theta<\pi$.
Moreover, the function
\begin{equation}
m(\theta) \equiv -\frac{1}{2}\ell'(\theta)= \csc^2\theta \,(1-\theta
\cot \theta) 
\label{eq:3.33}
\end{equation}
is positive, increasing, and convex for $0 < \theta < \pi$.  Also, 
$m(0)=1/3$, $m(\pi/2)=1$, and $m(\pi^-)=\infty$.
\end{lemma}

\begin{proof}
Using the product representation of $\sin \theta$, i.e., $\sin \theta
= \theta \prod_{k=1}^{\infty} (1-{\theta}^2/(k\pi)^2)$, one has 
\begin{equation}
\cot \theta = \sum_{k=-\infty}^{\infty} \frac{1}{\theta+k\pi}
\qquad
\mbox{and}
\qquad
\csc^2 \theta = \sum_{k=-\infty}^{\infty} \frac{1}{(\theta+k\pi)^2}
\label{eq:3.34}
\end{equation}
(convergence of the series for $\cot \theta$ here is understood in the 
sense of symmetric partial sums).  From (\ref{eq:3.34}) we obtain the 
following representation for $\ell(\theta)$
\begin{equation}
\ell(\theta)=- \sum_{k=1}^{\infty}\left[ \frac{k\pi}{(k\pi-\theta)^2}
- \frac{k\pi}{(k\pi+\theta)^2}\right].
\label{eq:3.35}
\end{equation}
It follows from (\ref{eq:3.35}) that $\ell(\theta) <0$ for
$0<\theta<\pi$.  Also from (\ref{eq:3.35}) we have 
\begin{equation}
m(\theta)=-\frac{\ell'(\theta)}{2}=
 \sum_{k=1}^{\infty}\left[ \frac{k\pi}{(k\pi+\theta)^3}
+ \frac{k\pi}{(k\pi-\theta)^3}\right],
\label{eq:3.36}
\end{equation}
which is positive for $0 <\theta <\pi$.  Thus, $m(\theta)$ is positive
and $\ell(\theta)$ is decreasing for $0<\theta<\pi$. Taking
derivatives again we find, 
\begin{equation}
m'(\theta)=
3 \sum_{k=1}^{\infty}\left[ \frac{k\pi}{(k\pi-\theta)^4}
- \frac{k\pi}{(k\pi+\theta)^4}\right].
\label{eq:3.37}
\end{equation}
The right-hand side of (\ref{eq:3.37}) is positive for $\theta \in
(0,\pi)$. Hence, $m(\theta)$ is increasing and $\ell(\theta)$ is
concave in $(0,\pi)$.  It also follows from (\ref{eq:3.37}) that
$m'(\theta)$ is increasing, and therefore $m(\theta)$ is convex in
$(0,\theta_1)$.  Lastly, the values of $m(\theta)$ at $\theta=0$, 
$\pi/2$, and $\pi^-$ are found by explicit evaluation.
\end{proof}

\begin{rem}
In fact the function $\ell$ and all its derivatives 
are negative for $0 < \theta < \pi$. 
\end{rem} 

\bigskip 

With all these preliminary results in hand we are ready to prove Theorems
\ref{thm:I2} and \ref{thm:I3}.

\bigskip

\begin{proof}[Proof of Theorem \ref{thm:I2}]
Since $\tilde \lambda_1(c)\equiv\lambda_1(\tilde
H_0(c))=c^2\lambda_1(c\theta_1)$ (see (\ref{eq:3.19}) above), we just
need to prove that   
\begin{equation}
\frac{d}{dc} \tilde \lambda_1(c) \Bigm|_{c=1} = \frac{1}{\theta_1}
\frac{d}{d\theta_1} \left(\theta_1^2 \lambda_1(\theta_1)\right) <0.  
\label{eq:3.38}
\end{equation}
This inequality follows from (\ref{eq:3.31}) and Lemmas \ref{lem:3.2}
and \ref{lem:3.3}.  Note that it holds for all $\theta_1 \in (0, \pi)$.
\end{proof} 

\bigskip

\begin{proof}[Proof of Theorem \ref{thm:I3}]
Since $\lambda_2(\theta_1) > \lambda_1(\theta_1)$ and $d \tilde
\lambda_1/dc\bigm|_{c=1} <0$, to prove (\ref{eq:3.25}) and therefore 
$$
\frac{d}{d\theta_1}
\left(\frac{\lambda_2(\theta_1)}{\lambda_1(\theta_1)}\right) >0
$$
reduces to showing (since ${\tilde \lambda}_1{\tilde \lambda}_2'- 
{\tilde \lambda}_2{\tilde \lambda}_1'$ can be grouped as ${\tilde
\lambda}_1 ({\tilde \lambda}_2'-{\tilde \lambda}_1')-({\tilde
\lambda}_2-{\tilde \lambda}_1) {\tilde \lambda}_1'$)
\begin{equation}
\frac{d\tilde \lambda_2}{dc}\Bigm|_{c=1} -  
\frac{d\tilde \lambda_1}{dc}\Bigm|_{c=1} >0, 
\label{eq:3.39}
\end{equation}
where $\tilde \lambda_2(c) \equiv \lambda_1(\tilde H_1(c))$. 
From equations (\ref{eq:3.28}) and (\ref{eq:3.29}) this is 
equivalent to proving
$$
(n+1) \frac{\int_0^{\theta_1} m(\theta) v_1^2 
\, d \theta}{\int_0^{\theta_1} v_1^2  \, d\theta} >
(n-3) \frac{\int_0^{\theta_1} m(\theta) v_0^2 \, d
\theta}{\int_0^{\theta_1} v_0^2  \, d\theta},
$$
which is obviously true for $n \le 3$ since $m(\theta)$ is positive.
For $n \ge 4$ it suffices to show
$$
\frac{\int_0^{\theta_1} m(\theta) v_1^2 
\, d \theta}{\int_0^{\theta_1} v_1^2  \, d\theta} >
\frac{\int_0^{\theta_1} m(\theta) v_0^2 \, d
\theta}{\int_0^{\theta_1} v_0^2  \, d\theta}.
$$
The fact that $g =v_1/v_0 = y_2/y_1$ is an increasing function of
$\theta$ on $[0,\theta_1]$ for $0 < \theta_1 \leq \pi/2$ (see 
Section 4 below) implies that $\hat v_1 \equiv
v_1/(\int_0^{\theta_1} v_1^2 \, d\theta)^{1/2}$ and $\hat v_0 \equiv
v_0/(\int_0^{\theta_1} v_0^2 \, d\theta)^{1/2}$ 
must have exactly one crossing in $(0,\theta_1)$.  Since $m(\theta)$
is positive and increasing in $[0,\pi)$ the desired inequality 
(\ref{eq:3.39}) follows by applying Lemma~2.7 of \cite{Ba80}, p.~69; 
this inequality is also known as Bank's inequality \cite{BaFl98} 
(see also \cite{Chi3} or \cite {AB92a}, p.~607).
\end{proof} 


\begin{rem}
For $n=2$ we can obtain a stronger version of Theorem \ref{thm:I3} 
which holds for all $\theta_1 \in (0,\pi)$.  Since $m(\theta)$ is 
increasing and $m(0)=1/3$ we have
$(n+1)m(\theta) - (n-1) \ge (n+1)/3-(n-1)=2(2-n)/3=0$ if $n=2$.  
Thus, from (\ref{eq:3.28}) we have 
$$
\frac{d \tilde \lambda_2}{dc}\Bigm|_{c=1} >0  \qquad \mbox{if $n=2$}.
$$
This inequality together with (\ref{eq:3.38}) implies
(\ref{eq:3.39}) (or, even more directly, (\ref{eq:3.25})) and 
therefore Theorem \ref{thm:I3} for $n=2$ and
$0<\theta_1<\pi$.  In fact, this argument shows that for all 
$\theta_1 \in (0,\pi)$ \, $\theta_1^2 \lambda_1(\theta_1)$ is 
decreasing and $\theta_1^2 \lambda_2(\theta_1)$ is increasing.  
These lead immediately to the fact that
$\lambda_2(\theta_1)/\lambda_1(\theta_1)$ and
$\theta_1^2[\lambda_2(\theta_1)-\lambda_1(\theta_1)]$ are both 
increasing for $0 < \theta_1 < \pi$.  For $n=3$ a related argument 
allows us to show that $\theta_1^2[\lambda_2(\theta_1)
-\lambda_1(\theta_1)]$ is increasing for $0 < \theta_1 < \pi$, 
and hence that Theorem \ref{thm:I3} extends to $0 < \theta_1 < \pi$ 
in that case as well.  For $n \ge 4$ we only have a proof of 
Theorem \ref{thm:I3} for $0 < \theta_1 \le \pi/2$. 
\end{rem}  

\bigskip

Next we prove two inequalities between the first two Dirichlet 
eigenvalues of a geodesic ball which are needed in Section 4.  


\begin{thm} \label{thm:3.1}
Let $\lambda_1$ and $\lambda_2$ be the first two eigenvalues of the
Dirichlet Laplacian on a geodesic ball contained in a hemisphere of
$\SS^n$. Denote its (geodesic) radius by $\theta_1$.  Then 
\begin{equation}
\frac{\lambda_2-n}{\lambda_1} \ge \frac{n+2}{n} \qquad 
\mbox{for $0 <\theta_1 \le \pi/2$}
\label{eq:3.40}
\end{equation}
with equality if and only if $\theta_1=\pi/2$ (i.e., for the hemisphere).
\end{thm} 

\begin{rems} 

(1) If we consider $\pi/2 < \theta_1 < \pi$ then inequality
(\ref{eq:3.40}) is reversed (and is a strict  inequality, i.e.,
$(\lambda_2 -n)/\lambda_1 < (n+2)/n$ for $\pi/2 < \theta_1 < \pi$).
This inequality is actually quite interesting since it allows us to
control the way $\lambda_2$ goes to $n$ (from above) as $\theta_1 \to
\pi^{-}$ in terms of $\lambda_1$ (since $\lambda_1$ goes to $0$ as
$\theta_1 \to \pi^{-}$).  The proof of this reversed inequality will
not be given, since it follows by making suitable modifications to
the proof of Theorem \ref{thm:3.1}, as given below.  

(2) Theorem \ref{thm:3.1} is the analog for $\SS^n$ of our Lemma 2.2 of
\cite{AB92b} for the Euclidean case.  Indeed, the inequality of Lemma
2.2 follows (except that we do not get a strict inequality) if we
consider (\ref{eq:3.40}) in the limit $\theta_1 \to 0^{+}$ (the
``Euclidean limit''). Since $\lambda_1 \approx \alpha^2/\theta_1^2$
and $\lambda_2 \approx \beta^2/\theta_1^2$ for $\theta_1$ near $0$
where $\alpha$ and $\beta$ are the Bessel function zeros
$j_{n/2-1,1}$ and $j_{n/2,1}$ (see \cite{AS64} for the notation here),
the inequality $\beta^2/\alpha^2 \ge (n+2)/n$ follows.    
\end{rems}  

\begin{proof}
The proof is similar to that of Lemma 2.2 of \cite{AB92b}.  We assume 
that $\theta_1 < \pi/2$ through most of the proof, returning to the 
case $\theta_1 = \pi/2$ (which can be treated explicitly in terms of 
elementary functions) only at the end.  We shall
use a suitable trial function (based on $y_2$) in the Rayleigh
quotient for $\lambda_1=\lambda_1(\theta_1)$: 
\begin{equation}
\lambda_1 \le \frac{\int_0^{\theta_1} (u')^2 \sin^{n-1} \theta \, d
\theta }{\int_0^{\theta_1} u^2 \sin^{n-1}\theta \, d \theta}
=-\frac{\int_0^{\theta_1} u( \sin^{n-1} \theta \, u')' \, d \theta
}{\int_0^{\theta_1} u^2 \sin^{n-1}\theta \, d \theta}.
\label{eq:3.41}
\end{equation}
It suffices that the trial function $u$ be real and continuous on
$[0,\theta_1]$, have $u(0)$ finite and $u(\theta_1)=0$, and be such
that all the integrals occurring above exist as finite real numbers.
This includes, in particular, the case $u=y_2/\sin \theta$, which we
now adopt. With this choice we find 
\begin{eqnarray}
-\left(\sin^{n-1}\theta \, u' \right)' &=& -\frac{d}{d\theta}
\left[\sin^{n-1} \theta \left(\frac{y_2}{\sin \theta} \right)' \right] \cr
 &=& -\sin^{n-2} \theta \, y_2''-(n-3) \sin^{n-2} \theta \cot \theta \, y_2'
\label{eq:3.42}\\
&+& (n-3) \sin^{n-4}\theta \, y_2 -(n-2) \sin^{n-2} \theta \, y_2, \cr 
\nonumber
\end{eqnarray}
and, upon using the differential equation satisfied by $y_2$ (i.e.,
(\ref{eq:3.1}) with $m=1$), we have 
\begin{equation}
-\left(\sin^{n-1}\theta \, u' \right)' =
2 \sin^{n-2} \theta \cot \theta \, y_2' +(\lambda_2-n+2) \sin^{n-2}
\theta \, y_2 - 2 \sin^{n-4} \theta \, y_2. 
\label{eq:3.43}
\end{equation}
That the expression on the right here has a finite limit as $\theta
\to 0^{+}$ follows from the fact that $y_2(\theta)$ can be taken as
$u_1(\theta)$ and $u_1(\theta)=c_1 \theta + O(\theta^3)$ as $\theta
\to 0^{+}$ (see (\ref{eq:3.1}) and (\ref{eq:3.2})).  It now follows from
(\ref{eq:3.41}) and (\ref{eq:3.43}) (since for $\theta_1 < \pi/2$,
$u=y_2/\sin \theta$ does not satisfy the equation that $y_1$  does, we
can write a strict inequality here) that 
\begin{equation}
\lambda_1 \int_0^{\theta_1} u y_2 \sin^{n-2} \theta \, d\theta <
\int_0^{\theta_1} u \left[2 \cot \theta \, y_2'+(\lambda_2-n+2) y_2 -2
\csc^2 \theta \, y_2 \right] \sin^{n-2} \theta \, d\theta. 
\label{eq:3.44}
\end{equation}
Since we would like to show that $\lambda_1 < n(\lambda_2-n)/(n+2)$,
it will be enough to show that the right-hand side of (\ref{eq:3.44})
is less than or equal to $[n(\lambda_2-n)/(n+2)] \int_0^{\theta_1} u
y_2 \sin^{n-2} \theta \, d \theta$, or, equivalently, 
$$
(n+2)\int_0^{\theta_1} u \left[2\cot \theta \, y_2' +(\lambda_2-n+2) y_2
-2 \csc^2 \theta \, y_2 \right] \sin^{n-2} \theta d \theta \le
n(\lambda_2-n) \int_0^{\theta_1} u y_2 \sin^{n-2} \theta d \theta 
$$
or 
\begin{equation}
2\int_0^{\theta_1} u \left[(n+2)\cot \theta \, y_2' +(\lambda_2+2) y_2
-(n+2) \csc^2 \theta \, y_2 \right] \sin^{n-2} \theta d \theta \le 0. 
\label{eq:3.45}
\end{equation}
We now rewrite this integral using $y_2=u_1(\theta;\lambda_2)$ (where
$\lambda_2=\lambda_2(\theta_1)$) and proceed to simplify the
integrand using the relations (\ref{eq:3.3}), (\ref{eq:3.6}) developed
for the $u_m$'s.  First, using (\ref{eq:3.3}) with $m=1$ the expression
in square brackets becomes 
$$
(n+2) \cot \theta [\cot \theta \, u_1 -u_2] + (\lambda_2+2) u_1 -(n+2)
\csc^2 \, \theta \, u_1 
$$
or, since $\cot^2  \theta  -\csc^2 \theta=-1$, 
$$
-(n+2) \cot \theta \, u_2 +(\lambda_2 -n) u_1.
$$
Finally, we use the recursion relation (\ref{eq:3.6}) with $m=2$ to 
see that
$$
u_3=(n+2) \cot \theta \, u_2 -(\lambda-n) u_1
$$
and hence, since we are taking $\lambda=\lambda_2$, inequality
(\ref{eq:3.45}) may be rewritten simply as 
$$
-2\int_0^{\theta_1} u u_3 \sin^{n-2} \theta \, d \theta \le 0,
$$
which clearly holds if $u_3>0$ for $0<\theta<\theta_1$ since
$u=y_2/\sin \theta$ and $y_2>0$ on $(0,\theta_1)$.  To see that
$u_3>0$ for  $0<\theta<\theta_1$ we use the relation (\ref{eq:3.7}),
first for $m=1$ and then for $m=2$.  With $m=1$ we have  
$$
\sin^n\theta \, u_2 (\theta)=(\lambda_2-n)\int_0^{\theta} (\sin t)^n
u_1(t) \, dt, 
$$
showing that $u_2(\theta;\lambda_2)>0$ for $0 < \theta \le \theta_1$
since $u_1(\theta;\lambda_2)=y_2(\theta) >0$ on $(0,\theta_1)$ and
$\lambda_2=\lambda_2(\theta_1)>n$ for $0<\theta_1 \le \pi/2$ by
domain monotonicity of Dirichlet eigenvalues (and the fact that
$\lambda_2(\pi/2)=2(n+1)>n$).  This  
in turn yields $u_3>0$ for $0<\theta\le \theta_1$ if
$0<\theta_1<\pi/2$ since by (\ref{eq:3.7}) with $m=2$ we have  
$$
\sin^{n+1}\theta \, u_3(\theta)=(\lambda_2-2(n+1))\int_0^{\theta} (\sin
t)^{n+1} u_2(t) \, dt, 
$$
and we know by domain monotonicity that $\lambda_2(\theta_1) >
\lambda_2(\pi/2)=2(n+1)$ for $0<\theta_1<\pi/2$.   

Thus we have proved inequality (\ref{eq:3.40}) of our theorem for
$0<\theta_1<\pi/2$ (with strict inequality in (\ref{eq:3.40})). 
It only remains to show
that equality holds when $\theta_1=\pi/2$ and this is elementary 
since in this case we can solve explicitly for all eigenvalues and
eigenfunctions without even the need of special functions.  One finds
that $\lambda_1=n$ with eigenfunction $x_{n+1}=\cos \theta$ and that
$\lambda_2=2(n+1)$ with exactly $n$ linearly independent
eigenfunctions $x_i x_{n+1}$ for $i=1,2, \dots, n$.  Clearly
$(\lambda_2-n)/\lambda_1=(n+2)/n$ and our proof is complete. 
\end{proof}  

\begin{rem}  
To fill in the picture for the hemisphere in $\SS^n$, we
note that $y_1=\cos \theta$ and $y_2=\sin \theta \cos \theta$, and
thus $u_1=\sin \theta \cos \theta$, $u_2=\sin^2 \theta$, and $u_3
\equiv 0$.  In this case $u=y_2/\sin \theta \equiv y_1$ and
(\ref{eq:3.41}) becomes an equality.  
\end{rem} 

\begin{lemma} \label{lem:3.4}
With notation as above,
\begin{equation}
\lambda_2 -\lambda_1 > \frac{n-1}{\sin^2 \theta_1}
\label{eq:3.46}
\end{equation}
for $0<\theta_1 \le \pi/2$.
\end{lemma}

\begin{proof} 
The proof is similar to that of Lemma 2.1 of \cite{AB92b}.  We use
$u=y_2$ in the Rayleigh quotient for $\lambda_1$ (i.e., in 
(\ref{eq:3.41}) above).  Since 
$$
-\left(\sin^{n-1} \theta \, y_2'\right)' = \sin^{n-1} \theta
\left(\lambda_2 -\frac{n-1}{\sin^2 \theta}\right) y_2,
$$
we get from (\ref{eq:3.41})
\begin{equation}
\lambda_1 < \lambda_2 - \frac{\int_0^{\theta_1} u^2
((n-1)/\sin^2\theta) \sin^{n-1} \theta \, d \theta}{\int_0^{\theta_1} 
u^2 \sin^{n-1} \theta \, d\theta} 
\label{eq:3.47}
\end{equation}
(again, $y_2$ does not satisfy the equation that $y_1$ does, so 
(\ref{eq:3.47}) is strict).  Now the lemma follows by the monotonicity
of $\sin \theta$ in $(0,\theta_1)$ for $0<\theta_1\le \pi/2$.  
\end{proof}   

To conclude this section we present the following results which are
needed in Section 4.   

\begin{thm} \label{thm:3.2}
Let
$p(\theta)=-u_0'(\theta,\lambda_1)/u_0(\theta,\lambda_1)=
u_1(\theta,\lambda_1)/u_0(\theta,\lambda_1)$.  Then, $p(\theta)$ is
positive, strictly increasing, and strictly convex on $(0,\theta_1)$
for $0<\theta_1\le \pi/2$.  Moreover, $p(0)=0$ and $p(\theta) \to
\infty$ as $\theta \to \theta_1^{-}$. 
\end{thm}

\begin{proof} The analog of this result in the Euclidean case was
proved in Lemma 2.3 of \cite{AB92b}.  That $p(0)=0$ and $p(\theta) \to
\infty$ as $\theta \to \theta_1^-$ follow from the boundary behavior of
$u_0$. Using the raising and lowering identities (\ref{eq:3.3}) and
(\ref{eq:3.4}) with $m=0$ and $m=1$ respectively and with $\lambda$ 
fixed at $\lambda_1(\theta_1)$, i.e., 
\begin{eqnarray}
u_1 &=& -u_0' \nonumber\\
& & \label{eq:3.48}\\
\lambda_1 u_0 &=& u_1'+(n-1) \cot \theta \, u_1, \nonumber
\end{eqnarray}
one obtains
\begin{equation}
u_0^2 \, p'=\lambda_1 u_0^2 +u_1^2 -(n-1)\cot \theta \, u_0 u_1 \equiv
\sigma(\theta), 
\label{eq:3.49}
\end{equation}
and in similar fashion
\begin{equation}
\left(\sin^{n-1} \theta \, \sigma(\theta)\right)'=(n-1)\sin^{n-3}
\theta \, u_0 u_1. 
\label{eq:3.50}
\end{equation}
The function $\sigma(\theta)$ is finite at $\theta=0$, and thus
$\sigma(\theta) \sin^{n-1} \theta=0$ at $\theta=0$. 
Since $u_0$ is positive and decreasing in $(0,\theta_1)$,
$u_0u_1=-u_0u_0'>0$ there.  Hence, $\sigma(\theta)>0$ in $(0,\theta_1]$
and therefore by (\ref{eq:3.49})  
$p$ is increasing there.  Clearly $p(\theta) > 0$ in $(0,\theta_1)$, since 
$p(0)=0$ and $p$ is increasing there.  Moreover, from (\ref{eq:3.49}) 
we find
\begin{equation}
u_0^2 p''=(n-1)\frac{u_0^2 p}{\sin^2 \theta}-(n-1) \cot \theta \, \sigma
+ 2 p \sigma \equiv s(\theta). 
\label{eq:3.51}
\end{equation}
From (\ref{eq:3.48}), (\ref{eq:3.49}), (\ref{eq:3.50}), and
(\ref{eq:3.51}) we obtain 
\begin{equation}
\left(\sin^{n+1} \theta \, s(\theta)\right)' = 2 \sin^n \theta\, \, 
\sigma(\theta) \, \Big{\{}\frac{\sigma}{u_0^2} \sin \theta+(n-1) \sin \theta 
+ 2 p \cos \theta \Big{\}}>0 
\label{eq:3.52}
\end{equation}
in $(0,\theta_1)$ for $0<\theta_1 \le \pi/2$. Since $\sin^{n+1}
\theta \,s(\theta)=0$ at $\theta=0$, this implies that
$s(\theta)>0$ in $(0,\theta_1)$ and therefore, by (\ref{eq:3.51}),
$p''>0$ there and the proof is complete. 
\end{proof} 

\begin{rem}
The function $p$ satisfies the Riccati equation
\begin{equation}
p'=\lambda_1+p^2-(n-1) \cot \theta \, p.
\label{eq:3.53}
\end{equation}
An alternative proof of Theorem \ref{thm:3.2} can be given directly from
(\ref{eq:3.53}).  In particular, the fact that $p$ is increasing follows
from (\ref{eq:3.53}) and the convexity of $\cot \theta$ in
$(0,\pi/2)$ using arguments similar to the ones used to prove Theorem 
\ref{thm:4.1} below.  
\end{rem} 

\bigskip

\begin{lemma} \label{lem:3.5}
The function $p(\theta)\cot \theta$ is strictly increasing on 
$(0,\theta_1)$ for $0<\theta_1 < \pi/2$. If $\theta_1=\pi/2$,
$p(\theta)=\tan\theta$, so $p(\theta)\cot \theta \equiv 1$ in that
case.   
\end{lemma}

\begin{proof}
Consider the function $r(\theta) \equiv p(\theta) \cot \theta
-\lambda_1/n$.  Using (\ref{eq:3.50}), (\ref{eq:3.51}), and
(\ref{eq:3.53}) one can show that $r(\theta)$ satisfies the equation 
\begin{equation}
\left(\sin^{n-1} \theta \, u_0^2 {} r'\right)'=2n \sin^{n-3} \theta
\, u_0^2 {} r. 
\label{eq:3.54}
\end{equation}
Since the function $p(\theta)$ is odd and analytic in $\theta$ for  
$\theta$ near $0$, we can expand it in odd powers of $\theta$. 
Inserting a series expansion in the Riccati equation (\ref{eq:3.53}) 
we find 
$$
p(\theta) = \frac{\lambda_1}{n} \theta +
\frac{\lambda_1}{3n^2(n+2)}\left(3\lambda_1 +n(n-1)\right)\theta^3
+O(\theta^5) 
$$
for $\theta$ near $0$. Substituting this into the definition of
$r(\theta)$ we find 
$$ 
r(\theta) \approx \frac{\lambda_1}{n^2(n+2)}(\lambda_1-n) \theta^2
$$
as $\theta \to 0^{+}$.  Since $\lambda_1(\theta_1) >n$ for $\theta_1 <
\pi/2$ (which follows from the fact that $\lambda_1(\pi/2)=n$ and
$\lambda_1(\theta_1)$ is decreasing in $\theta_1$), $r(\theta)$ is
positive in a neighborhood of $0$.  Now, $r(\theta)$ is a continuously
differentiable function in $[0,\theta_1)$, it is positive in a
neighborhood of $0$, and it goes to infinity as $\theta \to \theta_1^{-}$.
This implies that either $r(\theta)$ is an increasing function in
$[0,\theta_1)$ or that it has a positive local maximum.  However, this
latter possibility is ruled out by (\ref{eq:3.54}) (any possible
positive critical point of $r(\theta)$ must be a local minimum).
Thus, $r(\theta)$ is increasing on $(0,\theta_1)$ if $0<\theta_1<\pi/2$.  
If $\theta_1=\pi/2$, $u_0=\cos \theta$, hence $u_0'=-\sin\theta$,
$p(\theta)=\tan \theta$, and finally $r(\theta)\equiv 0$ since 
$\lambda_1=n$.  
\end{proof}  

\begin{rem}
In the Euclidean case, the analog of Lemma \ref{lem:3.5}
(i.e., Lemma 2.4 of \cite{AB92b}) was an immediate consequence of
the convexity of the function analogous to the function $p(\theta)$
used here.  For $\SS^n$, as we have seen, the proof is somewhat more
involved.  
\end{rem}  

\bigskip

\section{Monotonicity properties of $g$ and $B$}  

\bigskip

    In this section we prove the key properties of the functions that
occur in our rearrangement procedure in Section 6.  These functions are
\begin{equation}
q(\theta) = \sin \theta \frac{g'(\theta)}{g(\theta)},
\label{eq:4.1} 
\end{equation}
and
\begin{equation}
B(\theta)= g'(\theta)^2+\frac{n-1}{\sin^2 \theta}
g(\theta)^2=[q^2+(n-1)]\left(\frac{g}{\sin \theta}\right)^2,
\label{eq:4.2}
\end{equation}
where 
\begin{equation}
g(\theta)=\frac{y_2(\theta)}{y_1(\theta)}=
\frac{u_1(\theta;\lambda_2(\theta_1))}{u_0(\theta;\lambda_1(\theta_1))}.
\label{eq:4.3}
\end{equation}
Our objectives in this section are to prove that $g(\theta)$ is increasing
and $B(\theta)$ is decreasing in $[0,\theta_1]$ for $0 < \theta_1 \le
\pi/2$.  For the hemisphere (i.e., for $\theta_1=\pi/2$) we can
explicitly compute  $g(\theta)=\sin \theta$ and
$B(\theta)=(n-1)+\cos^2 \theta$ (see the remark following the proof
of  Theorem \ref{thm:3.1}).  It is obvious that $g$ is increasing 
and $B$ is decreasing in this case.  Thus, we can assume in the rest
of this section that $\theta_1 < \pi/2$. 

Since $g$ and $\sin \theta$ are positive, that $g$ is increasing will
be a simple consequence of showing $q \ge 0$.  On the other hand, from
(\ref{eq:4.1}) and (\ref{eq:4.2}) it follows that
\begin{equation}
B'(\theta)=2 \left [ qq'-(\cos \theta -q)(q^2+n-1) \frac{1}{\sin
\theta} \right] \left(\frac{g(\theta)}{\sin \theta} \right)^2.
\label{eq:4.4} 
\end{equation}
Hence, that $B$ is decreasing will be a consequence of showing that
$q'\le 0$ and $0 \le q \le \cos \theta$.  Thus, in order to prove the
desired properties of $g$ and $B$ we only need to show that $0 \le
q(\theta) \le \cos \theta$ and $q'(\theta) \le 0$ for $0 \le \theta
\le \theta_1$.  

The strategy we use to prove these results for $q$
follows the same general method used in \cite{AB92b} and \cite{AB95}.  
Since we shall need the boundary behavior of $q$ at the two endpoints
$\theta=0$ and $\theta=\theta_1$ (this is necessitated by the fact 
that the coefficients in the right-hand side of the differential equation 
for $q$, equation (\ref{eq:4.7}) below, become singular at the two 
endpoints), we give these now.  By Taylor--Frobenius expansion, we find
\begin{equation}
q(0)=1, \qquad  q'(0)=0, \qquad q''(0)=2
\left(\frac{\lambda_1}{n}-\frac{\lambda_2}{n+2} -\frac{2-n}{2(n+2)}
\right), 
\label{eq:4.5}
\end{equation}
\begin{equation}
q(\theta_1)=0, \qquad  q'(\theta_1)=-\frac{1}{3}
\left(\lambda_2-\lambda_1- \frac{n-1}{\sin^2 \theta_1}\right) \sin
\theta_1. 
\label{eq:4.6}
\end{equation}
Throughout this section we will use $\lambda_1$ and $\lambda_2$ to denote 
$\lambda_1(\theta_1)$ and $\lambda_2(\theta_1)$, respectively.  From 
(\ref{eq:4.5}) and Theorem \ref{thm:3.1} we have that $q''(0) < -1$ for 
$\theta_1 < \pi/2$.  From (\ref{eq:4.6}) and Lemma \ref{lem:3.4} we have 
that $q'(\theta_1) < 0$.  Therefore, $q< \cos \theta$ on an interval just 
to the right of $0$ and $q>0$ on an interval just to the left of $\theta_1$.  

In order to prove that $q \ge 0$, $q' \le 0$, and $q \le \cos \theta$
for $0 \le \theta \le \theta_1$ we analyze the ordinary differential
equation satisfied by $q$.  To obtain the differential equation for
$q$ first we differentiate (\ref{eq:4.1}) with respect to $\theta$,
and use our choice of $g=y_2/y_1$.  Then we use the equation
(\ref{eq:3.1}), with $m=0$ and $\lambda=\lambda_1(\theta_1)$,
satisfied by $y_1$, and the same equation, but this time with $m=1$ and
$\lambda=\lambda_2(\theta_1)$, satisfied by $y_2$.  Thus, we obtain 
\begin{equation}
q'=2pq-(n-2) q \cot \theta - \frac{q^2+1-n}{\sin \theta} -
(\lambda_2 -\lambda_1) \sin \theta.
\label{eq:4.7} 
\end{equation}
Here, the function $p \equiv - y_1'/y_1$ obeys the Riccati equation
\begin{equation}
p'-p^2 +(n-1) \cot \theta \,  p - \lambda_1=0
\label{eq:4.8} 
\end{equation}
associated to equation (\ref{eq:3.1}) with $m=0$ and 
$\lambda=\lambda_1(\theta_1)$ (note that (\ref{eq:4.7}) is also a
Riccati equation).  We proved in Section 3 that $p$ is positive,
strictly increasing, and strictly convex in $(0,\theta_1)$, for 
$0 < \theta_1 \le \pi/2$, $p(0)=0$, and $p \to \infty$ as $\theta 
\to \theta_1^{-}$ 
(see Theorem \ref{thm:3.2}). 

Having derived the equation for $q$, we are ready to prove the
necessary facts about $q$.  We start by showing $q \ge 0$ in
$[0,\theta_1]$, which we prove by contradiction. Assume $q$ is
negative somewhere in $[0,\theta_1]$.  Since $q(0)=1$ and $q$ is
positive to the left of $\theta_1$ (and $q$ is continuous) this
implies that there are two points $\alpha$, $\beta$, say, with
$0<\alpha<\beta<\theta_1$ such that $q(\alpha)=q(\beta)=0$ and
$q'(\alpha) \le 0$, $q'(\beta) \ge 0$. At points in $(0,\theta_1)$
where $q=0$ it follows from (\ref{eq:4.7}) that
\begin{equation}
q'=\frac{n-1}{\sin \theta} - (\lambda_2 - \lambda_1) \sin \theta.
\label{eq:4.9} 
\end{equation}
Since $\lambda_2 > \lambda_1$ and $\sin \theta$ is increasing in
$\theta$ for $0<\theta<\pi/2$, the right-hand side of (\ref{eq:4.9}) is 
strictly decreasing in $\theta$, and hence it is not possible to have 
$\alpha <\beta$ with $q(\alpha)=q(\beta)=0$ and $q'(\alpha) \le 0$, 
$q'(\beta) \ge 0$. Therefore $q \ge 0$ in $[0,\theta_1]$. 

The proof that $q \le \cos \theta$ follows the same ideas. 
Define the function $\psi = \cos \theta -q$. From (\ref{eq:4.7}) we get
\begin{equation}
\psi'=(\lambda_2-\lambda_1-n) \sin \theta + \Big( 2 \, p +\frac{\psi}{\sin
\theta} \Big) (\psi-\cos \theta)-(n-1)\psi \cot \theta.
\label{eq:4.10}
\end{equation}
We have already shown that $\psi(\theta)=\cos \theta-q(\theta)$ is
positive in a neighborhood of $\theta=0$.  Also, if $\theta_1 <\pi/2$,
$\psi(\theta_1)=\cos \theta_1 >0$. Now assume that $\psi$ is negative
somewhere in $[0,\theta_1]$.  This implies that there are two points
$r$, $s$, say, with $0<r<s<\theta_1$ such that $\psi(r)=\psi(s)=0$ and
$\psi'(r) \le 0$, $\psi'(s) \ge 0$. At points in $(0, \theta_1)$
where $\psi=0$ we have from (\ref{eq:4.10}) that
\begin{equation}
\frac{1}{\sin \theta} \psi'=(\lambda_2-\lambda_1-n) -2 \, p \, 
\cot \theta.
\label{eq:4.11}
\end{equation}
From Lemma \ref{lem:3.5} it follows that the right-hand side of
(\ref{eq:4.11}) is strictly decreasing in $\theta$, hence it is not
possible to have $r<s$ with $\psi(r)=\psi(s)=0$ and $\psi'(r) \le 0$,
$\psi'(s) \ge 0$. Thus $\psi \ge 0$, and hence $q \le \cos \theta$ in
$[0,\theta_1]$.    

Finally, we show that $q(\theta)$ is decreasing in $[0,\theta_1]$. 
It is convenient now to write (\ref{eq:4.7}) in an alternative form
and use convexity arguments.  Specifically we consider
\begin{equation}
q'=2\, p(\theta) \, q +(n-2)(1-q) \cot \theta +\frac{1-q^2}{\sin \theta}
+ (n-2) \tan(\theta/2)- (\lambda_2-\lambda_1) \sin \theta
\label{eq:4.12}
\end{equation}
and, since $\lambda_2-\lambda_1>0$, $q \le 1$, and $n \ge 2$, the
right-hand side of (\ref{eq:4.12}), $F(\theta,q)$, say, is convex in
$\theta$ for fixed $q$ because each of the functions $\cot \theta$,
$\csc \theta$, $\tan(\theta/2)$, and $-\sin \theta$ is individually
convex on the interval $[0,\theta_1]$ for $0 < \theta_1 \le \pi/2$.
Also the function $p(\theta)$ is convex on the same interval (see
Theorem \ref{thm:3.2} above).  To see how these facts imply that $q'
\le 0$ observe that if not we could find three points $\alpha_1$,
$\alpha_2$, $\alpha_3$ in $(0,\theta_1)$, where
$q(\alpha_1)=q(\alpha_2)=q(\alpha_3)$, $q'(\alpha_1) <0$,
$q'(\alpha_2)>0$, and $q'(\alpha_3) <0$. But then we would have (we use
$q$ for the common value of $q(\alpha_i)$ here)
\begin{eqnarray}
0 < q'(\alpha_2) &=& F(\alpha_2,q)=F(\mu \, \alpha_1 + (1-\mu) \alpha_3,q)
\nonumber\\
&<& \mu F(\alpha_1,q)+(1-\mu)F(\alpha_3,q)= \mu
q'(\alpha_1)+(1-\mu)q'(\alpha_3) <0,
\label{eq:4.13}
\end{eqnarray}
a contradiction.  We have used the strict convexity of $F(\theta,q)$
in $\theta$ here.  The parameter $\mu$ is strictly between $0$ and $1$
and determines $\alpha_2$ as a convex combination of $\alpha_1$ and
$\alpha_3$, that is,  $\alpha_2=\mu \, \alpha_1+(1-\mu) \alpha_3$. 
To summarize this section we state the results above as a theorem. 

\begin{thm} \label{thm:4.1} 
With $q(\theta)$, $B(\theta)$, and $g(\theta)$ as defined in
(\ref{eq:4.1}), (\ref{eq:4.2}), and (\ref{eq:4.3}), respectively, the
inequalities $0 \le q \le \cos \theta$ and $q' \le 0$ hold on
$[0,\theta_1]$ for $0 < \theta_1 \le \pi/2$.  It follows that
$g(\theta)$ is increasing and $B(\theta)$ is decreasing there as well.
\end{thm}

\bigskip

\section{Chiti's comparison argument in $\SS^n$} 

\bigskip 

    Here we need an extension of Chiti's comparison result \cite{Chi0}, 
\cite{Chi1}, \cite{Chi2}, \cite{Chi3} (see also Appendix A of \cite{AB92a}), 
given originally for domains in $\RR^n$, to the case of domains in 
$\SS^n$.  

We let $\SS^n(\rho)$ denote the $n$--dimensional sphere of radius
$\rho$ (hence of constant sectional curvature $\kappa=1/\rho^2$).
$\SS^n$ will always denote $\SS^n(1)$.  Define the function 
$$
S_{\rho}(r)=\rho \sin(r/\rho).
$$
It is well known that, in geodesic polar coordinates, the metric on
$\SS^n(\rho)$ is
$$
ds^2=dr^2+S_{\rho}(r)^2 {\vert d \omega \vert}^2,
$$
where $r$ represents geodesic distance from a point and ${\vert d \omega 
\vert}^2$ is the canonical metric for $\SS^{n-1}$.  One should think of
$r$ as $\rho$ times the angle $\theta$ from the north pole.  Thus
$r$ runs from $0$ to $\rho \pi$. 

We will always assume that $r=0$ (the north pole) is the center for our 
spherical rearrangements (this is not a restriction since the metric is
identical in geodesic polar coordinates about any point).  Then for a
bounded domain $\Omega \subset \SS^n(\rho)$ {} $\Omega^{\star}$, the
spherical rearrangement of $\Omega$, will denote the geodesic ball
about $r=0$ having the same volume as $\Omega$, i.e., $\vert
\Omega^{\star} \vert=\vert \Omega \vert$.   
 
For a nonnegative function $f$ defined on $\Omega$ we define two 
rearranged functions $f^{\#}$ and $f^{\star}$.  The decreasing 
rearrangement $f^{\#}$ of $f$ is a function from $[0,\vert \Omega 
\vert]$ to $\RR$ which is equimeasurable with $f$ and nonincreasing. 
We will use $s$ as the argument of $f^{\#}$ in most instances.  The 
symmetric decreasing (or spherical decreasing) rearrangement 
$f^{\star}$ of $f$ is a function defined on $\Omega^{\star}$ which 
is invariant under rotations about $r=0$, equimeasurable with $f$, 
and nonincreasing with respect to $r$. $f^{\star}$ is a function of 
$x \in \Omega^{\star}$ but because of its symmetry we will abuse 
notation and write $f^{\star}(r)$ where $r$ is the geodesic distance 
from the center of $\Omega^{\star}$.  With this understanding we have 
$$
f^{\star}(r)=f^{\#}(A(r))
$$
where
\begin{equation}
A(r) \equiv s = n C_n \int_0^r S_\rho(\tau)^{n-1} \, d \tau
\label{eq:5.1} 
\end{equation}
is the $n$-volume of the geodesic ball of radius $r$ in
$\SS^n(\rho)$. Here $C_n=\pi^{n/2}/\Gamma(\frac{n}{2}+1)$ denotes the
volume of the unit ball in $\RR^n$ (and $n C_n$ is its ``surface area'', 
i.e., $|\SS^{n-1}|$).  We shall also have occasion to
use increasing rearrangements.  These will be denoted  
$f_{\#}$ (the increasing rearrangement of $f$) and $f_{\star}$ (the 
symmetric or spherical increasing rearrangement of $f$) and their 
definitions are analogous to those of $f^{\#}$ and $f^{\star}$, 
respectively.  For further information on rearrangements the reader 
is referred to the book of Hardy, Littlewood, and P\'{o}lya \cite{HLP52} 
and the many other references given in \cite{AB92a},\cite {Ta76b}. 
  
Finally, we also need the classical isoperimetric inequality extended 
to $\SS^n(\rho)$ in its sharp form.  We let $L(t)$ be the function
giving the $(n-1)$--dimensional volume of the geodesic ball of radius 
$r$, i.e.,
\begin{equation}
L(r)=n C_n S_{\rho}(r)^{n-1}=A'(r).     
\label{eq:5.2} 
\end{equation}
Then, for example, for any domain in $\SS^2$ one has
\begin{equation}
L^2\ge 4\pi A -(A/\rho)^2 = 4 \pi A - \kappa A^2      
\label{eq:5.3} 
\end{equation}
(see, e.g., \cite{Os78}).  Here $A=\vert \Omega \vert$ and $L$ is the 
length of the boundary of $\Omega$.  Equality occurs in (\ref{eq:5.3}) 
if and only if $\Omega$  is a geodesic ball.  In $\SS^n$, $n>2$, the 
sharp classical isoperimetric  inequality cannot be given in as explicit 
a form as (\ref{eq:5.3}).  For a bounded domain $\Omega$ we define 
$H_{n-1}(\partial \Omega)$ to be the $(n-1)$--dimensional volume of 
$\partial \Omega$.  The isoperimetric inequality on $\SS^n(\rho)$ then 
reads
\begin{equation}
H_{n-1}(\partial \Omega) \ge H_{n-1}(\partial \Omega^{\star})
\label{eq:5.4} 
\end{equation}
with equality if and only if $\Omega$ is a geodesic ball (see Burago 
and Zalgaller \cite{BuZa88}, p.~86, Theorem 10.2.1).  In terms of the 
function $L(r)$ defined above (\ref{eq:5.4}) may be written 
\begin{equation}
H_{n-1}(\partial \Omega) \ge L(\theta(\vert\Omega \vert))= n C_n 
S_{\rho}(\theta(\vert \Omega \vert))^{n-1},    
\label{eq:5.5} 
\end{equation}
where $\theta(s)$ is the inverse function to the function $A$ defined
in (\ref{eq:5.1}). 

With all these ingredients we state our extension of Chiti's 
comparison result (throughout the rest of this section we set
$\rho=1$, the extension to arbitrary $\rho$ being straightforward).
\begin{thm} \label{thm:5.1}
Let $\Omega$ be a bounded domain in $\SS^n$ and let $\lambda_1$ and
$u_1$ denote the first Dirichlet eigenvalue and eigenfunction of the 
Laplacian on $\Omega$.  Let $B_{\lambda_1}$ be the geodesic ball of such 
a radius that $\lambda_1$ is also the first Dirichlet eigenvalue of 
the Laplacian on $B_{\lambda_1}$.  Let $v_1 > 0$ be the first Dirichlet  
eigenfunction on $B_{\lambda_1}$ and fix its normalization so that 
$\int_{\Omega} u_1^2 \, d\sigma= \int_{B_{\lambda_1}} v_1^2 \, d\sigma$.  

Then there is a value $r_1 \in (0,\theta(\vert B_{\lambda_1} \vert))$ 
such that
\begin{equation}
v_1(r) \ge u_1^{\star}(r)  \qquad \mbox{for $r\in [0,r_1]$},   
\label{eq:5.6}
\end{equation}
and
\begin{equation}
v_1(r) \le u_1^{\star}(r)  \qquad \mbox{for $r\in [r_1,\theta(\vert 
B_{\lambda_1} \vert)]$},    
\label{eq:5.7} 
\end{equation}
where $\theta(s)$ is the function defined following (\ref{eq:5.5}) above. 
\end{thm} 

\begin{proof}
Let $u_1$, respectively $\lambda_1$, be the lowest
eigenfunction, respectively eigenvalue, of the Dirichlet problem on
$\Omega \subset \SS^n$, i.e.,
\begin{equation}
-\Delta u_1=\lambda_1 u_1  \qquad \mbox{in $\Omega$},
\qquad u_1=0  \qquad \mbox{on $\partial \Omega$.}    
\label{eq:5.8} 
\end{equation}
Define $\Omega_t=\{x \bigm| u_1(x)>t \}$ and $\partial \Omega_t=\{x
\bigm| u_1(x)=t \}$. Let $\mu_1(t)=\vert \Omega_t\vert$, and $\vert
\partial \Omega_t \vert\equiv H_{n-1}(\partial \Omega_t)$, where $H_{n-1} 
(d \sigma)$ denotes $(n-1)$--dimensional measure on $\SS^n$.  Then we
have (see, e.g., Talenti \cite{Ta76b}, p.~709, eq.~(32))
\begin{equation}
-\mu_1'(t)=\int_{\Omega_t} \frac{1}{\vert \nabla u_1 \vert} H_{n-1}
(d\sigma),
\label{eq:5.9} 
\end{equation}
for almost every $t>0$.  Applying Gauss's theorem to (\ref{eq:5.8}), 
we have
\begin{equation}
\lambda_1 \int_{\Omega_t} u_1 \, d\sigma= \int_{\partial \Omega_t}
\vert \nabla u_1 \vert H_{n-1}(d\sigma),
\label{eq:5.10} 
\end{equation}
since the outward normal to $\Omega_t$ is given by $-\nabla u_1/|\nabla u_1|$. 
Using the Cauchy--Schwarz inequality and equations (\ref{eq:5.9}) and
(\ref{eq:5.10}) we obtain 
\begin{equation}
{\vert \partial \Omega_t \vert}^2 = {\left( \int_{\partial \Omega_t}
H_{n-1}(d \sigma) \right)}^2 \le (-\mu_1'(t)) \lambda_1
\int_{\Omega_t} u_1 \, d \sigma.
\label{eq:5.11} 
\end{equation}
As discussed above, if $\Omega$ is a domain in $\SS^n$, the
classical isoperimetric inequality is given by
\begin{equation}
H_{n-1} (\partial \Omega) \ge H_{n-1} (\partial \Omega^{\star}) 
\label{eq:5.12} 
\end{equation}
where $\Omega^{\star}$ is a geodesic ball having the same volume as
$\Omega$. The $(n-1)$--dimensional measure of $\partial \Omega^{\star}$,
$H_{n-1}(\partial \Omega^{\star})$, is given in terms of
$\theta_1(\Omega^{\star})$, the geodesic radius of $\Omega^{\star}$, by
\begin{equation}
H_{n-1}(\partial \Omega^{\star})=nC_n {\left(\sin \theta_1(\Omega^{\star})
\right)}^{n-1}, 
\label{eq:5.13} 
\end{equation}
where $C_n=\pi^{n/2}/\Gamma(\frac{n}{2}+1)$ is the volume of the unit
ball in $\RR^n$ (and $nC_n$ is the volume of $\SS^{n-1}$).
Therefore, from (\ref{eq:5.12}) and (\ref{eq:5.13}) it follows that
\begin{equation}
\vert \partial \Omega_t \vert =  H_{n-1}(\partial \Omega_t) \ge nC_n
{\left(\sin \theta_1(\Omega_t^{\star}) \right)}^{n-1}. 
\label{eq:5.14} 
\end{equation}
Hence, from (\ref{eq:5.11}) we have
\begin{equation}
\lambda_1 \int_{\Omega_t} u_1 \, d \sigma \ge n^2 C_n^2
{\left[ \sin \theta_1 (\Omega_t^{\star}) \right]}^{2n-2} \left( -
\frac{1}{\mu_1'(t)} \right). 
\label{eq:5.15} 
\end{equation}
Finally one uses the fact that
\begin{equation}
\int_{\Omega_t} u_1 \, d \sigma= \int_0^{\mu_1(t)} u_1^{\#} (s) \, ds,
\label{eq:5.16} 
\end{equation}
which follows directly from the definition of $u_1^{\#}$, the decreasing
rearrangement of $u_1$ on the interval $[0, \vert \Omega \vert]$
(here $s$ is a variable denoting volume and is related to the
geodesic radial variable $\theta$ via $s=nC_n \int_0^{\theta} (\sin
r)^{n-1} \, dr$, i.e., $ds/d\theta=nC_n(\sin \theta)^{n-1}$). 

Since $u_1^{\#}(s)$ is the inverse function to $\mu_1 (t)$ we have
$$
-\frac{d u_1^{\#}}{ds} = - \frac{1}{\mu_1'(t)},
$$
which, combined with (\ref{eq:5.15}) and (\ref{eq:5.16}), yields
\begin{equation}
-\frac{d u_1^{\#}}{ds} \le \lambda_1 n^{-2}C_n^{-2} (\sin
\theta(s))^{2-2n} \int_0^s u_1^{\#}(s') \, ds'.
\label{eq:5.17} 
\end{equation}
Now, if we view $v_1$ as a function of the volume $s$ (here we will abuse
notation and continue to call it $v_1$) rather than as a function of
$\theta$ (or $r$), with $s=n C_n \int_0^{\theta} (\sin r)^{n-1} \, dr$, it satisfies
(\ref{eq:5.17}) with equality, i.e.,
\begin{equation}
-\frac{d v_1}{ds} = \lambda_1 n^{-2}C_n^{-2} (\sin
\theta(s))^{2-2n} \int_0^s v_1 (s') \, ds'.
\label{eq:5.18} 
\end{equation}
In fact, (\ref{eq:5.18}) is an integrated version of equation
(\ref{eq:3.1}) with $m=0$ and $\lambda=\lambda_1$ in the variable $s$.
Having obtained the relations (\ref{eq:5.17}) and (\ref{eq:5.18})
for $u_1^{\#}(s)$ and $v_1(s)$ respectively, we will prove that under
the normalization condition imposed on them, they are either
identical or they cross only once in the interval $(0,\vert
B_{\lambda_1} \vert)$ (in the sense specified by Theorem \ref{thm:5.1}). 
All the arguments we give below depend on the
continuity of $u_1^{\#}$ and $v_1$.  The function $v_1$ is in fact real
analytic in $[0,\vert B_{\lambda_1} \vert]$ and, furthermore, it is 
decreasing there as can be seen from Lemma \ref{lem:3.2} above 
(or from (\ref{eq:5.18}) and the fact that $v_1>0$).  The absolute continuity
of $u_1^{\#}$ on $[0,\vert \Omega \vert]$ follows from arguments in 
\cite{Ta76b}.  Since $u_1^{\#}$ and $v_1$ are normalized to have the same
$L^2$--norm, they either are identical or they cross.  If they
coincide, then $B_{\lambda_1}=\Omega^{\star}$, and Theorem \ref{thm:5.1} 
is proved since any $r \in (0,\theta(\vert B_{\lambda_1} \vert))$ will 
serve as $r_1$. Next, 
following Chiti \cite{Chi2} we conclude that $u_1^{\#}(0)$ cannot exceed 
$v_1(0)$.  In fact, if $u_1^{\#}(0) \ge v_1(0)$ it follows by mimicking 
the proof of the main theorem in  \cite{Chi2} that $v_1(s) \le u_1^{\#}(s)$ 
for all $s \in [0,\vert B_{\lambda_1} \vert)$ which, in turn, implies 
$v_1(s) \equiv u_1^{\#}(s)$ and we are back in the previous case.  Thus, if
$v_1(s) \not \equiv  u_1^{\#}(s)$, $v_1(s) > u_1^{\#}(s)$ in a neighborhood of
$0$, and both functions being of the same norm they must cross at
least once.  Let $s_1$ be the largest $s \in (0, \vert B_{\lambda_1}
\vert)$ such that  $u_1^{\#}(s') \le v_1(s')$ for all $s'\le s$.  By the
definition of $s_1$, there is an interval immediately to the right of
$s_1$ on which $u_1^{\#}(s) > v_1(s)$.  Indeed, by continuity and the
definition of $s_1$ 
\begin{equation}
\int_{0}^{s} [u_1^{\#}(s')-v_1(s')] \, ds' < 0
\qquad \mbox{for $0 < s \le s_1 + \epsilon$}
\label{eq:5.18a}
\end{equation}
for some $\epsilon>0$.  It now follows that $u_1^{\#}(s) > v_1(s)$ at least
on the interval from $s_1$ to $s_1+\epsilon$ since by the absolute
continuity of $u_1^{\#}$
\begin{eqnarray}
v_1(s)-u_1^{\#}(s) &=& \int_{s_1}^{s} \left[\frac{d}{ds}(v_1-u_1^{\#})
\right] \, ds
\label{eq:5.18b} \\
&\le& \lambda_1 n^{-2}C_n^{-2} \int_{s_1}^{s} {\left(\sin \theta(s')
\right)}^{2-2n} \int_0^{s'} [u_1^{\#}(s'')-v_1(s'')] \, ds'' \, ds'
\nonumber\\
&<&0 \qquad \mbox{for $s \in (s_1,s_1+\epsilon]$} \nonumber
\end{eqnarray}
by virtue of (\ref{eq:5.18a}).

We will now show that $u_1^{\#}(s) > v_1(s)$ for all $s \in (s_1, 
\vert B_{\lambda_1} \vert]$, which will prove the theorem.  If not, 
then the point $s_2$ defined to be the largest $s \in (s_1, 
\vert B_{\lambda_1}\vert]$ for which $u_1^{\#}(s') > v_1(s')$ for 
all $s_1 < s' < s$, would be less than $\vert B_{\lambda_1} \vert$, 
and we would have
\begin{equation}
u_1^{\#}(s) > v_1(s)  \qquad \mbox{for $s\in (s_1,s_2)$}
\label{eq:5.19}
\end{equation}
with $u_1^{\#}(s_1)= v_1(s_1)$ and $u_1^{\#}(s_2)= v_1(s_2)$.  In this 
case, we can define the function
\begin{equation}
w(s) = \begin{cases}
v_1(s)  & \mbox{for $s \in [0,s_1] \cup [s_2, \vert B_{\lambda_1}
\vert]$} \\
u_1^{\#}(s) & \mbox{for $s \in (s_1,s_2)$.}
\end{cases}
\nonumber
\end{equation}
It follows from (\ref{eq:5.17}) and (\ref{eq:5.18}) that $w$ satisfies
\begin{equation}
-\frac{d w}{ds}(s)  \le  \lambda_1 n^{-2}C_n^{-2} (\sin
\theta(s))^{2-2n} \int_0^s w (s') \, ds',
\label{eq:5.20} 
\end{equation}
with strict inequality for all $s>s_1$.  From $w(s)$ define the function 
$$g(x)=w(s(\theta))$$
for $x \in B_{\lambda_1}$ where $\theta$ is the polar angle (angle from 
the north pole) corresponding to $x$.  Thus $g$ is a radial function on
$B_{\lambda_1}$ (assumed centered at the north pole).  Because of 
(\ref{eq:5.20}) (or (\ref{eq:5.18b}) or (\ref{eq:5.19})), $g$ cannot 
be the groundstate of the Laplacian with Dirichlet boundary conditions 
on $B_{\lambda_1}$ (but it is certainly an admissible trial function 
for $\lambda_1$).  Therefore, 
\begin{equation}
\lambda_1 < \frac{\int_{B_{\lambda_1}} \vert \nabla g \vert^2 \,
d\sigma}{\int_{B_{\lambda_1}} g^2 \, d\sigma}.
\label{eq:5.21}
\end{equation}
By standard change of variables,
\begin{equation}
\int_{B_{\lambda_1}} g^2 \, d\sigma= \int_0^{\vert B_{\lambda_1} \vert}
w^2(s) \, ds
\label{eq:5.22}
\end{equation}
and
\begin{equation}
\int_{B_{\lambda_1}} \vert \nabla g \vert^2 \, d\sigma= 
n^2 C_n^2 \int_0^{\vert B_{\lambda_1} \vert} (\sin \theta(s))^{2n-2} 
w'(s)^2 \, ds.
\label{eq:5.23}
\end{equation}
Using (\ref{eq:5.20}) (substitute for one of the $w'(s)$'s in
(\ref{eq:5.23}), using the fact that $-w'(s)>0$) and integration by
parts we get 
\begin{equation}
\int_{B_{\lambda_1}} \vert \nabla g \vert^2 \, d\sigma \le 
\lambda_1 \int_0^{\vert B_{\lambda_1} \vert}
w(s)^2 \, ds.
\label{eq:5.24}
\end{equation}
From (\ref{eq:5.21}), (\ref{eq:5.22}), and (\ref{eq:5.24}) we get a
contradiction, and the theorem follows.
\end{proof} 

\bigskip

\section{The main result} 
\bigskip 

    After all the preliminaries developed in Sections 2 through 5 we are 
ready to prove our main result, i.e., Theorem \ref{thm:I1} from 
which the PPW result for domains contained in a hemisphere of $\SS^n$ 
follows as indicated in the introduction.  As in our proof of the PPW 
conjecture for domains in $\RR^n$, the starting point here is the use 
of the {\it gap inequality}, which is a variational estimate for the 
difference between the first two eigenvalues of the Laplacian.  The gap 
inequality states that
\begin{equation}
\lambda_2(\Omega) -\lambda_1 (\Omega) \le 
\frac{\int_{\Omega} {\vert \nabla P \vert}^2
u_1^2 \, d \sigma}
{\int_{\Omega} P^2 u_1^2 \, d \sigma}   
\label{eq:6.1} 
\end{equation}
provided $\int_{\Omega} P u_1^2 \, d \sigma=0$ and $P \not \equiv 0$.
Here $\Omega$ is a domain in $\SS^n$ and $d \sigma$ is the standard 
volume element in $\SS^n$. 
The gap inequality follows from  the Rayleigh--Ritz inequality for
$\lambda_2$ using $Pu_1$ as the trial function (hence the side
condition $Pu_1 \perp u_1$) after a suitable integration by parts.  To
get the desired isoperimetric result out of this one must make very
special choices of the function $P$: in particular, choices such that
(\ref{eq:6.1}) is an equality if $\Omega$ is the appropriate geodesic 
ball.  


A key element needed to guarantee the orthogonality of the trial 
functions $P_i u_1$ to $u_1$ which we will use in the sequel 
is the {\it center of mass} argument embodied in Theorem
\ref{thm:CM1} above.


Concerning the choice of trial functions $P_i$ we proceed as follows. 
Thinking of $\SS^n$ as the unit sphere in $\RR^{n+1}$ and with the 
{\it center of mass} point for $\Omega$ fixed at the north pole we take
\begin{equation}
P_i=g(\theta) \frac{x_i}{\sin \theta},  \qquad \mbox{$i=1,2,\dots,n$},
\label{eq:6.2}
\end{equation}
where $\theta$ represents the angle of a point from the positive
$x_{n+1}$--axis (the direction of the north pole). 
In $\SS^n$ the variable $\theta$ is the geodesic
radial variable with respect to the north pole.  Division by $\sin
\theta=\sqrt{1-x_{n+1}^2}$ normalizes the $n$--vector $(x_1,x_2,\dots,x_n)$. 
The choice of $g(\theta)$, as in the Euclidean case, is determined by
the fact that we must have equality in (\ref{eq:6.1}) when $\Omega$ is a
geodesic ball.  For a geodesic ball of radius $\theta_1$, $u_1=c_1y_1$
where $y_1$ satisfies 
\begin{equation}
{y_1}''+(n-1) \cot \theta \,  y_1' + \lambda_1 y_1 =0
\label{eq:6.3}
\end{equation}
with boundary conditions $y_1(0)$ finite and $y_1(\theta_1)=0$. This 
is just equation (\ref{eq:3.1}) with $m=0$ and $\lambda=\lambda_1$.  Also, 
$u_2=c_2 (x_i/\sin \theta) y_2$ (for any $i=1,2,\dots,n$) where $y_2$
satisfies
\begin{equation}
y_2''+(n-1) \cot \theta \, y_2' +\Big(\lambda_2-\frac{n-1}{\sin^2 \theta}\Big)
y_2 =0
\label{eq:6.4} 
\end{equation}
with boundary conditions $y_2(0)=y_2(\theta_1)=0$ (this is just equation
(\ref{eq:3.1}) with $m=1$ and $\lambda=\lambda_2$).  The values
$\lambda_1$ and $\lambda_2$ are to be taken as the least eigenvalues
of these one-dimensional radial problems.  By Lemma \ref{lem:3.1} these 
are the correct identifications of $\lambda_1$ and $\lambda_2$ for our 
geodesic ball.  One can express $y_1$ and $y_2$ in terms of associated 
Legendre functions.  In fact,
$$
y_1(\theta) = (\sin \theta)^{1-n/2} \, P_{\nu_1}^{-(n/2-1)}(\cos \theta)
$$
and
$$
y_2(\theta) = (\sin \theta)^{1-n/2} \, P_{\nu_2}^{-n/2}(\cos \theta)
$$
up to unimportant constant factors, where the parameters $\nu_1$ and 
$\nu_2$ are related to the eigenvalues $\lambda_1$ and $\lambda_2$ 
respectively by
$$
\nu(\nu+1)=\lambda+\frac{n(n-2)}{4}.
$$ 
We follow Abramowitz and Stegun \cite {AS64} in our notation here; note 
that their convention for associated Legendre functions makes 
$(\sin \theta)^{1-n/2} \, P_{\nu}^{-\mu}$, where $\mu=n/2-1+m$ and $m$ 
is a nonnegative integer, and not $(\sin \theta)^{1-n/2} \, P_{\nu}^\mu$, 
the ``right'' $n$-dimensional generalization of the familiar associated 
Legendre functions $P_{\nu}^m$, $m=0,1,\ldots$, from $\SS^2$ (or $\RR^3$).  
This means, in particular, that when $n$ is even $P_{\nu}^{\mu}$ can be 
substituted for $P_{\nu}^{-\mu}$ (they are then proportional), while if 
$n$ is odd $Q_{\nu}^{\mu}$ can be used (the distinction here is between 
$\mu$ being an integer or half an odd integer).  Using $P_{\nu}^{-\mu}$ 
with $\mu$ as above circumvents this ``even-odd effect''.  Since we want 
equality in (\ref{eq:6.1}) when $\Omega$ is a geodesic ball, using the 
form of $u_1$ and $u_2$ for a geodesic ball we see that $g(\theta)$ must 
be essentially the quotient of $y_2$ by $y_1$. 

Let $\Omega \subset \SS^n$ be contained in a hemisphere.  Let 
$B_{\lambda_1}$ denote the geodesic ball in $\SS^n$ having the same 
value of $\lambda_1$ as $\Omega$.  By Sperner's inequality \cite{Sp73} 
(see also \cite{FrHa76}) $\lambda_1(B_{\lambda_1})=\lambda_1(\Omega)\ge
\lambda_1(\Omega^{\star})$ and by the properties of $\lambda_1$ for
geodesic balls (in particular, $\lambda_1$ decreases as the radius of
the ball increases), we see that $\theta_1 \le \pi/2$, where
$\theta_1$ denotes the geodesic radius of $B_{\lambda_1}$. 
We now set $g(\theta)=y_2(\theta)/y_1(\theta)$ for $0 \le \theta \le
\theta_1$, $g(\theta) = g(\theta_1)$ for $\theta_1 \le \theta \le
\pi/2$ and we extend $g(\theta)$ to a function on $[0,\pi]$ by
reflection about $\theta=\pi/2$.  (Note that this definition makes 
$B(\theta) = g'(\theta)^2 + (n-1) g(\theta)^2 / \sin^2 \theta$ a 
decreasing function on $[\theta_1, \pi/2]$, since $B(\theta) 
= (n-1) g(\theta_1)^2/\sin^2 \theta$ there.) We then apply the 
{\it center of mass} result (Theorem \ref{thm:CM1}) to $\Omega$
and $\tilde G (\theta)=g(\theta)/\sin \theta$ obtaining a choice of
Cartesian coordinate axes such that
$$
\int_{\Omega} g(\theta)\frac{x_i}{\sin \theta} u_1^2 \, d\sigma=0, \qquad 
\mbox{for $i=1,2, \dots, n$.}
$$
Using the functions $P_i=g(\theta) x_i/\sin \theta$ in (\ref{eq:6.1}) we get
\begin{equation}
(\lambda_2 - \lambda_1) \int_{\Omega} g(\theta)^2
{\left( \frac{x_i}{\sin \theta}\right)}^2 u_1^2 \, d\sigma \le
\int_{\Omega} {\Bigl\vert  
\nabla \left( g(\theta) \frac{x_i}{\sin \theta}\right) \Bigr\vert}^2 u_1^2
\, d\sigma
\label{eq:6.5}
\end{equation}
for $1 \le i \le n$, and summing on $i$ from $1$ to $n$ we obtain
\begin{equation}
\lambda_2-\lambda_1 \le \frac{\int_{\Omega} B(\theta) u_1^2 \, d\sigma}
{\int_{\Omega} g(\theta)^2 u_1^2 \, d\sigma},
\label{eq:6.6}
\end{equation}
where
\begin{equation}
B(\theta) = {g'(\theta)}^2+\frac{n-1}{\sin^2 \theta} g(\theta)^2
\label{eq:6.7}
\end{equation}
as defined previously (see equation (\ref{eq:4.2})).  We observe that, 
in spite of the fact that $\Omega$ is contained in a hemisphere, our 
use of the center of mass result may imply that $\Omega$ does not lie 
in the northern hemisphere (i.e., $\theta$ would not be in $(0,\pi/2)$ 
for all points in $\Omega$).  To remedy this situation, we observe that 
since $\Omega$ is contained in a hemisphere, $-\Omega \subset \SS^n
\setminus \Omega$. Thus, if we let $\Omega_{\pm}= \{ \vec x \in
\Omega \bigm| \pm x_{n+1} >0 \}$ we have $\Omega_{+} \cap (-\Omega_{-}) 
= \emptyset$. Since $g(\theta)$ and $B(\theta)$ are both symmetric
with respect to $\theta= \pi/2$ it follows that the integrals over
$\Omega$ can be replaced by integrals over  $\tilde \Omega = \Omega_{+} 
\cup  (-\Omega_{-} )$ with no change in their values if we agree to
transplant $u_1$ along with $\Omega_{-}$ to $-\Omega_{-}$. This
follows since $g$ and $B$ were defined to be symmetric about
$\theta=\pi/2$ and therefore they transplant into themselves.
After moving the whole problem to the northern hemisphere (where $g$
is increasing and $B$ is decreasing) we can carry out all the further
rearrangements in exact parallel with the Euclidean case,
encountering no further difficulties. 

To conclude the proof of Theorem \ref{thm:I1} we need the following
two chains of inequalities.  We have
\begin{eqnarray}
\int_{\Omega} B(\theta) u_1^2 \, d\sigma =\int_{\tilde {\Omega}} B(\theta) 
\tilde u_1^2 \,d\sigma 
&\leq& \int_{{\Omega}^{\star}} 
B(\theta)^{\star}u_1^{\star 2} \,d\sigma \nonumber \\
&\leq& \int_{{\Omega}^{\star}} B(\theta)
u_1^{\star 2}\, d\sigma \leq \int_{B_{\lambda_1}} B(\theta) v_1^2 \, d\sigma,
\label{eq:6.8}
\end{eqnarray}
and
\begin{eqnarray}
\int_{\Omega} g(\theta)^2 u_1^2 \, d\sigma = \int_{\tilde {\Omega}} 
g(\theta)^2 \tilde u_1^2 \, d\sigma 
&\geq& \int_{{\Omega}^{\star}} g(\theta)^2_\star u_1^{\star 2} \, d \sigma
\nonumber \\
&\geq& \int_{{\Omega}^{\star}} g(\theta)^2 
u_1^{\star 2} \, d \sigma \geq \int_{B_{\lambda_1}} g(\theta)^2 v_1^2 \, 
d\sigma.
\label{eq:6.9}
\end{eqnarray}
Here $\tilde u_1$ represents $u_1$ as transplanted to $\tilde {\Omega}$ 
and $v_1$ is the first eigenfunction of $-\Delta$ on $B_{\lambda_1}$ 
with Dirichlet boundary conditions and normalized so that 
$\int_{\Omega} u_1^2 \, d \sigma =\int_{B_{\lambda_1}} v_1^2 \, d\sigma$. 
The functions $g$ and $B$ are likewise based on the eigenfunctions of the 
ball $B_{\lambda_1}$ (so that the $\theta_1$ that goes into the boundary 
value problems (\ref{eq:6.3}) and (\ref{eq:6.4}) that define them is the 
radius of the ball $B_{\lambda_1}$).  In each of (\ref{eq:6.8}) and 
(\ref{eq:6.9}), the equality is trivial, the first inequality follows 
simply from rearrangement (see Section 5 for our notation), the
second inequality is by virtue of the monotonicity properties of 
$g(\theta)$ and $B(\theta)$, and the last inequality follows from our 
$\SS^n$ analog of Chiti's comparison result (see Section 5 for details) 
and also uses the monotonicity properties of $g$ and $B$ again.  Finally, 
from (\ref{eq:6.6}), (\ref{eq:6.8}), and (\ref{eq:6.9}) we obtain
\begin{equation}
\lambda_2(\Omega) -\lambda_1(\Omega) \leq \frac{\int_{B_{\lambda_1}} 
B(\theta) v_1^2 \, d\sigma}
{\int_{B_{\lambda_1}}g(\theta)^2 v_1^2 \, d\sigma} =\lambda_2
(B_{\lambda_1})- \lambda_1(B_{\lambda_1}).
\label{eq:6.10} 
\end{equation}
Hence (since $\lambda_1(B_{\lambda_1})=\lambda_1(\Omega)$)
\begin{equation}
\lambda_2(\Omega) \leq \lambda_2 (B_{\lambda_1})
\label{eq:6.11}
\end{equation}
which concludes the proof of Theorem \ref{thm:I1},
it being clear from any of a number of our previous inequalities that 
equality obtains in (\ref{eq:6.11}) if and only if $\Omega$ is itself 
a ball. 

\bigskip

\begin{rems}

(1) Theorems \ref{thm:I1} (i.e., inequality (\ref{eq:6.11})) 
and \ref{thm:I4} also hold under somewhat more general 
circumstances than for $\Omega$ contained in a hemisphere of 
$\SS^n$. In particular, they continue to hold if $\Omega \cap 
(-\Omega)=\emptyset$, or, more generally, if $\Omega$ satisfies the 
``excess less than or equal to deficit property'' with respect 
to the center of mass as north (or south!)\ pole as stated in 
(\ref{eq:CM12}) above for all $k \in [0,1]$.
We note, however, that even under these conditions $\Omega$ is 
constrained to have volume no larger than that of a hemisphere. 
See also Remark 4 in Section 2 following the proof of Theorem 
\ref{thm:CM1}.

(2) In fact, it is enough that $\Omega$ satisfy $\theta_1(B_{\lambda_1})
\le \pi/2$ (or, equivalently, $\lambda_1(\Omega) \ge n$), where
$\theta_1(B_{\lambda_1})$ is the geodesic radius of the ball $B_{\lambda_1}$,
together with the ``excess less than or equal to deficit property''
(\ref{eq:CM12}) for all $k \in [\cos \theta_1(B_{\lambda_1}),1]$.  
This condition allows us to prove that $\lambda_2(\Omega) \le 
\lambda_2(B_{\lambda_1})$ even for certain domains which have volume 
larger than that of a hemisphere (as well as covering all previous 
cases).  A variant of this condition also applies in the case of the 
Neumann problem (the ``$\mu_1$ problem'') for domains in $\SS^n$ (see 
\cite{AB95} and certain of our remarks in Section 2 above).  Then 
$\theta_1(\Omega^\star)$, the geodesic radius of $\Omega^\star$, should 
replace $\theta_1(B_{\lambda_1})$ in the foregoing (and in this case 
we are still limited by $|\Omega| \le \frac{1}{2}|\SS^n|$).  The reason 
for these values of $\theta_1$ is that these are the radii we end with 
in the respective problems, when all is said and done.    
\end{rems}


\section*{Acknowledgements}

\bigskip

   M.S.A. is grateful for the hospitality of Thomas Hoffmann--Ostenhof
and the Erwin Schr\"{o}\-dinger Institute (ESI) in Vienna, where some of
this work was carried out.  We would also like to thank the referee for 
several useful remarks.

\end{document}